\documentclass[10pt]{article}
\usepackage{amsfonts, epsfig, amsmath, amssymb, color}

\usepackage[english]{babel}

\newcommand{\I}[1]{\mathbb{I}\{#1\} }
\newcommand{\E}{\mathbf{E}}
\renewcommand{\P}{\mathbf{P}}
\newcommand{\Px}{\mathbf{P}_{\!x}}

\renewcommand {\epsilon}{\varepsilon}

\newtheorem{theorem}{Theorem}%[section]
\newtheorem{defin}{Definition}[section]
\newtheorem{prop}{Proposition}[section]
\newtheorem{cor}{Corollary}[section]

\newtheorem{lemma}{Lemma}[section]

\newtheorem{rem}{Remark}[section]

\DeclareMathSymbol{\ophi}{\mathalpha}{letters}{"1E}

\newcommand{\g}{\gamma}
\newcommand{\e}{\varepsilon}
\renewcommand{\phi}{\varphi}

\newcommand{\be}{\begin{equation}}
\newcommand{\ee}{\end{equation}}
\newcommand{\ben}{\begin{equation*}}
\newcommand{\een}{\end{equation*}}

\newcommand{\ba}{\begin{equation}\begin{aligned}}
\newcommand{\ea}{\end{aligned}\end{equation}}

\DeclareMathOperator{\Const}{Const}

\DeclareMathOperator{\sgn}{sgn}

\DeclareMathOperator{\Law}{Law}

\newenvironment{proof}{\par\noindent{\bf Proof:}}{\hfill$\blacksquare$\par}
\newenvironment{Aproof}[2]{\par\noindent{\bf Proof of #1 #2:}}{\hfill$\blacksquare$\par}

\newfont{\cyrfnt}{wncyr10}
\def\J3{\cyrfnt{\rm \u{\cyrfnt I}}}
\def\j3{\cyrfnt{\rm \u{\cyrfnt i}}}

\allowdisplaybreaks[4]

\begin{document}
%\layout
\title{Metastable Behaviour of Small Noise L\'evy-Driven Diffusions
}

%\date{\today\footnote{{\jobname}.tex\hfill \textbf{Preliminary version}}}

\date{\today}

\author{Peter Imkeller and Ilya Pavlyukevich}

\maketitle

\begin{abstract}
We consider a dynamical system in $\mathbb{R}$ driven by a vector field $-U'$, where $U$ is a multi-well
potential satisfying some regularity conditions. We perturb this dynamical system by a L\'evy noise
of small intensity and such that the heaviest tail of its L\'evy measure is regularly varying. 
We show that the perturbed dynamical system exhibits metastable behaviour i.e.\ on a proper time scale
it reminds of a Markov jump process taking values in the local minima of the potential $U$. Due to the heavy-tail
nature of the random perturbation, the results differ strongly from the well studied purely Gaussian case.  
\end{abstract}



\tableofcontents
\numberwithin{equation}{section}

\section{Introduction}

This paper addresses the rigorous mathematical 
description of the phenomenon of metastability in systems 
with big jumps. The picture we shall study may be outlined as follows. Let us consider a one-dimensional
deterministic dynamical system driven by a vector field $-U'(\cdot)$, where $U(\cdot)$ is a multi-well 
potential with some smoothness conditions and a certain increase rate at infinity. 
According to the initial conditions the deterministic trajectories of the dynamical 
system converge to the local minima of the potential $U$ or stay in its local maxima. Obviously, no transition
between different domains of attraction is possible.

The situation becomes different if the dynamical system is perturbed by (small) random noise whose presence  allows transitions between the potential wells. However depending on the system's initial conditions and 
noise's properties, 
certain potential wells may be reached only on appropriately long time scales or stay unvisited. 
The phenomenon of metastability means, roughly speaking, that for different time scales 
and initial conditions the system may reach different local statistical equilibria. 

The system's behaviour is determined by the type of random perturbation.
Unquestionably, dynamical systems subject to small Gaussian perturbations have been studied most extensively.
The main reference on this subject is the book \cite{FreidlinW-98} where the large deviations theory for the
perturbed trajectories is established.   
The large deviations estimates allow to solve the first exit problem from 
the domain of attraction of a stable point. It turns out that the mean exit time is exponentially large
in the small noise parameter, and its logarithmic rate is proportional to the height of the potential 
barrier the trajectories have to overcome. Thus for a multi-well dynamical system we obtain 
a series of exponentially non-equivalent time scales given by the wells' mean exit times. Moreover, one
can prove that the normalised exit times are exponentially distributed (see \cite{Williams-82, Day-83, BovierEGK-04}), and thus
have a memoryless property which is referred to in physical literature as unpredictability.

In the simplest situation when the potential $U$ has only two wells of different depths, 
one can observe two statistically different regimes.
First, if the time horizon is shorter than the exit time from the shallow well, the system cannot leave the well 
where it has started, and therefore stays in the neigbourhood of 
the well's local minimum. Second, if the time horizon is longer than the exit time from the shallow well, 
the system has enough time to reach the deepest well from any starting point, and stays in the vicinity 
of the global minimum. In \cite{KipnisN-85} the following metastability result is established. Namely, there is
a time scale on which the dynamical system converges to a Markov two-state process with one absorbing state 
corresponding to the deep well. It is easy to notice that this particular time scale is given by the mean 
exit time from the shallow well. More general results for multidimensional diffusions can be found in \cite{Mathieu-95} and \cite{GalvesOV-87}.  

There is a very close connection between metastability of a small noise system and spectral properties of its 
infinitesimal operator. It can be shown that exponentially small eigenvalues of the infinitesimal
generator
are expressed in terms of mean life times in the domains of attraction, and the corresponding
eigenfunctions are close to constants on these domains \cite{KolokoltsovM-96}. On the other hand, the generator's
eigenvalues can be calculated with the help of variational principles 
\cite{BuslovM-92, BovierGK-05}.

However, recently non-Gaussian perturbations with big jumps attract more attention. 
Instant transitions between remote states are referred to as extreme events and are observed in
dynamics of asset prices, climate and telecommunication systems etc.
In the physical literature, non-Gaussian symmetric stable L\'evy processes are
used especially often, under the name of \textit{L\'evy flights}. 
The mathematical study of the gradient dynamical systems subject to small perturbation by a heavy-tail L\'evy
process was tackled in \cite{ImkellerP-06} (for symmetric stable processes), 
where results on the first exit time form the potential
well with non-characteristic boundary were established by purely probabilistic methods. 
It was shown that the exit time increases as a power of the small noise 
parameter and does not depend on the
depth of the potential well but rather on the distance between the local minimum and the domain's boundary.

In the present paper which can be seen as a sequel of \cite{ImkellerP-06} we deal with more general multi-well
potential and arbitrary L\'evy processes with regularly varying tails. 
The presence of big jumps makes the L\'evy driven dynamics quite different from the purely
Gaussian one. Indeed, the life times in the potential wells belong now to the same time scale which leads to a quite
different process in the limit of small parameter.

%\newpage

\section{Object of study and main result}

Let $(\Omega,\mathcal{F},(\mathcal{F})_{t\geq 0},\P)$ be a filtered probability space.
We assume that the filtration satisfies the usual hypotheses in the sense of
\cite{Protter-04}, i.e.\ $\mathcal{F}_0$ contains all the 
$\P$-null sets of $\mathcal{F}$, and is right continuous.

We consider solutions $X^\e=(X_t^\e)_{t\geq 0}$ of the one-dimensional stochastic differential equation
\be
\label{eq:main}
X^\e_t(x)=x-\int_0^t U'(X_{s-}^\e(x))\, ds+\e L_t,\quad x\in \mathbb{R},
\ee
where $L$ is a L\'evy process and $U$ is a potential function satisfying the following assumptions.

\medskip
Assumptions on $L$:
\begin{enumerate}
\item[\textbf{L1}] 
$L$ has a generating triplet $(d,\nu,\mu)$ with a Gaussian variance $d\geq 0$, an arbitrary drift $\mu\in\mathbb{R}$
and a L\'evy measure $\nu$ satisfying the usual condition $\int_{\mathbb{R}\backslash\{0\}}\max\{y^2,1\}\,\nu(dy)<\infty$.
For $u\geq 1$ denote the tails of the L\'evy measure $\nu$
\ba
H_-(-u)=\int_{(-\infty,-u)}\nu(dy),\qquad
H_+(u)=\int_{(u,+\infty)}\nu(dy),
\ea
and $H(u)=H_-(-u)+H_+(u)$.
\item[\textbf{L2}] 
Assume, $H_+(\cdot)$ is regularly varying at infinity, i.e.\
\be
H_+(u)=u^{-r}l(u), \quad u\to+\infty,
\ee
for some $r>0$ and a slowly varying function $l$ (for regular variation see Appendix \ref{a:rv}).
\item[\textbf{L3}] 
Assume that there exists a finite limit
\ba
\lim_{u\to+\infty}\frac{H_-(-u)}{H_+(u)}=\kappa\in (0,+\infty).
\ea
or
\ba
\limsup_{u\to+\infty}\frac{H_-(-u)}{H_+(u)}=\kappa=0.
\ea
\end{enumerate}
\medskip
\noindent
Assumptions on $U$:
\begin{enumerate}
\item[\textbf{U1}] 
$U\in\mathcal{C}^1(\mathbb{R})\cap\mathcal{C}^3([-K,K])$
for some $K>0$ large enough.
\item[\textbf{U2}]  
$U$ has exactly $n$ local minima $m_i$, $1\leq i\leq n$, and $n-1$ local maxima $s_i$, 
$1\leq i\leq n-1$, enumerated in increasing order 
\be
-\infty=s_0<m_1<s_1<m_2<\cdots<s_{n-1}<m_n<s_n=+\infty.
\ee
All extrema of $U$ are non-degenerate, i.e.\ 
$U''(m_i)>0$, $1\leq i\leq n$, and $U''(s_i)<0$, 
$1\leq i\leq n-1$.
\item[\textbf{U3}]  
$|U'(x)|> c_1 |x|^{1+c_2}$ as $x\to\pm\infty$ for
some $c_1$, $c_2>0$.
\end{enumerate}

The class of L\'evy processes $L$ under consideration covers for example 
compound Poisson processes with heavy-tail jumps or stable L\'evy processes with L\'evy measure
\be
\nu(dy)=\left( c_1\I{y<0}+c_2\I{y>0}\right) \frac{dy}{|y|^{1+\alpha}},\quad \alpha\in(0,2),
\quad c_1\geq 0, c_2>0.
\ee
%
% It is easy to see, that if $\kappa >0$, the negative tail $H_-(\cdot)$ is regularly varying at
%$-\infty$ with the same
% index $-r$.

We consider $X^\e$ for small values of $\e$, $\e\downarrow 0$.

Since the L{\'e}vy process $L$ is a semimartingale, the stochastic differential equation  
\eqref{eq:main} is well defined, see also \cite{Protter-04} for the general theory.
However, since the drift term $U'$ is not globally Lipschitz we need to show the existence 
and uniqueness of the strong solution of \eqref{eq:main} which is done in Appendix \ref{a:ss}.
 
Under assumptions on $U$, the underlying deterministic ($\e=0$) equation
\be
\label{eq:main0}
X^0_t(x)=x-\int_0^t U'(X_s^0(x))\, ds
\ee
has a unique solution for any initial value $x\in \mathbb{R}$ and all $t\geq 0$. 
The local minima of $U$ are stable attractors for the dynamical system $X^0$, i.e.\
if $x\in(s_{i-1},s_i)$, 
$1\leq i\leq n$, then $X^0_t(x)\to m_i$ as $t\to\infty$. It is clear that the deterministic solution $X^0$
does not leave the domain of attraction where it started.

Our goal is to describe the phenomenon of \textit{metastability} which roughly speaking consists in
the existence
of a time scale for which the system reminds of a jump process taking values in the set stable attractors.  
We prove the following main Theorem.

\begin{theorem}
\label{th:main}
Let $X^\e(x)=(X^\e_t(x))_{t\geq 0}$ be a solution of \eqref{eq:main}. If $x\in (s_{i-1}, s_i)$,
for some $i=1,\dots,n$, then
\ba
X^\e_{t/ H(1/\e)}(x)\to Y_t(m_i),\quad \e\downarrow 0,
\ea
in the sense of finite-dimensional distributions, where $Y=(Y_t)_{t\geq 0}$ is a Markov process on a state
space $\{m_1,\dots, m_n\}$ with the infinitesimal generator $Q=(q_{ij})_{i,j=1}^n$,
\ba
\label{eq:q}
q_{ij}&
=
\frac{\kappa\I{j<i}+\I{j>i}}{1+\kappa}\left||s_{j-1}-m_i|^{-r}-|s_j-m_i|^{-r}\right|,\quad i\neq j,\\
-q_{ii}&=q_i=\sum_{j\neq i}q_{ij}=\frac{\kappa}{1+\kappa}|s_{i-1}-m_i|^{-r}
+\frac{1}{1+\kappa}|s_i-m_i|^{-r} .
\ea

\end{theorem}

Let us consider a particular example of equation \eqref{eq:main}, namely 
a symmetric $\alpha$-stable process $L$ (L\'evy flights) in a double-well
potential. Let $U$ satisfy Assumptions formulated above and let for definiteness 
$s_1=0$. The process $L$ has 
a generating triplet $(0,\nu,0)$ with a L\'evy measure $\nu(dy)=|y|^{-1-\alpha}$, $y\neq 0$,
$\alpha\in (0,2)$. Such dynamics is often considered in physical literature.  
P.\ Ditlevsen in \cite{Ditlevsen-99a,Ditlevsen-99b} studied such a system
in his attempt to explain abrupt catastrophic climate changes during the last Ice Age.
Further in \cite{ChechkinGKM-05}, the authors addressed the calculation of the mean transition
time between the wells if $\alpha\in [1,2)$. 
(Their conclusions based on numerical simulations of the process $X^\e$ 
are not fully consistent with our results, and thus should be improved.)
One-well dynamics of such processes was firstly studied in \cite{ImkellerP-06}.

Due to Theorem \ref{th:main}, the main
features of the process $X^\varepsilon$ in the small noise limit are
retained by a Markov jump process, and on the time scale
$\alpha\varepsilon^{-\alpha}$ we obtain the following convergence in the
sense of finite dimensional distributions:
\begin{equation}
X^\varepsilon_{\alpha t/\e^{\alpha}}(x)\to Y_t, \quad t>0,\quad\varepsilon\downarrow 0,
\end{equation}
where $Y$ is a Markov process on the state space $\{m_1,m_2\}$ with the
following matrix as infinitesimal generator
\begin{equation}
Q=
\begin{pmatrix}
 -m_1^{-\alpha}& \hphantom{-}m_1^{-\alpha} \\
\hphantom{-}m_2^{-\alpha} & -m_2^{-\alpha}
\end{pmatrix}
\quad{\rm and }
\quad
Y_0=
\begin{cases}
m_1, & \text{ if } x<0,\\ m_2, & \text{ if } x>0.
\end{cases}
\end{equation}

To compare the result obtained with its Gaussian
counterpart we refer to \cite{KipnisN-85}, where this problem
was first studied.

Let us consider a Gaussian diffusion $\hat{X}^\varepsilon$
which solves the equation
\be 
\label{eq:xg}
\hat{X}^\e_t(x)=x-\int_0^t U'(\hat{X}^\e_s(x))\, ds+\e W_t,
\ee
where $W$ is a standard Brownian motion.
Since it is well known that in
the Gaussian case the height of the potential barriers plays a
crucial role, we assume that the left well is deeper, i.e.\ $U(0)-U(m_1)>U(0)-U(m_2)$.
This leads to the following
meta-stable behaviour of $\hat{X}^\varepsilon$ (\cite[Theorem 2.1]{KipnisN-85}). These exists a time scale
$\lambda^\e$ such that 
\be
\lim_{\e\to 0} \e^2\ln\lambda^\e=2(U(0)-U(m_2))
\ee
and
\begin{equation}
\hat{X}^\varepsilon_{t\lambda^\varepsilon}(x)\to
\hat{Y}_t,
\quad\varepsilon\downarrow 0,
\end{equation}
in the sense of finite dimensional distributions, where
$\hat{Y}$ is a Markov process on $\{m_1,m_2\}$
with the infinitesimal matrix
\begin{equation}
\begin{pmatrix}
 0& \hphantom{-}0 \\
1 & -1
\end{pmatrix}
\quad{\rm and }\quad
\hat{Y}_0=
\begin{cases}
m_1, & \text{ if } x<0,\\
m_2, & \text{ if } x>0.
\end{cases}
\end{equation}
As we see, the main difference between L\'evy and Gaussian dynamics
consists not only in different intrinsic time scales --- polynomial
vs.\ exponential, --- but also in a qualitatively different limiting
behaviour. In the heavy-tail case, the states of the limiting process
are recurrent, whereas in the Gaussian case, the minimum of the
deepest well is absorbing.

In general case, we can summarise the differences as follows.
First, we see that the characteristic time scale
is algebraic in $\e$. 
Second, the properties of the limiting process $Y$ depend on sizes of the potential wells
and not on their depths. Further, if $\kappa>0$, the all states of $Y$ are recurrent. 
The process $Y$ has a unique absorbing state $m_n$ 
(the local minimum of the right peripheral well) if and only if $\kappa=0$, i.e.\ 
when the positive tail of $L$ dominates.

% Let us briefly compare our results with a purely Gaussian case. First, we see that the characteristical time scale
% is subexponential in $\e$. Second, the properties of the limiting process $Y$ depend on sizes of the potential wells
% and not on thier depths. Further, if $\kappa>0$, the all states of $Y$ are recurrent. The process $Y$ has 
% a unique absorbing state $m_n$ 
% (the local minimum of the right peripheral well) if and only if $\kappa=0$, i.e.\ 
% when positive tail of $L$ dominates.  

\medskip

This material is organised as follows. In Section \ref{s:onewell} we decompose the L\'evy process $L$ into
a small jump part and a compound Poisson part and study the small-jump dynamics of the process $X^\e$. 
Section \ref{s:law:sigma} is devoted to the asymptotics of the
first exit time from a single well. Section \ref{s:transitions} provides the asymptotic exponentiality of the
transition times between the wells. Theorem~\ref{th:main} is proved in Section \ref{s:meta}. 
Appendices \ref{a:ss} and \ref{a:rv} contain the proof of the existence of the strong solution 
of \eqref{eq:main} and basic information on regularly varying functions. 

\medskip

\textsc{Acknowledgments:} 
This work was supported by the DFG Research Project `Stochastic Dynamics of Climate States',
DFG Research Center MATHEON in Berlin and the Japan Society for the Promotion of Science.
I.\ Pavlyukevich thanks Th.\ Mikosch, G.\ Samorodnitsky, Th.\ Simon and 
N.\ Yoshida for their valuable comments.

%%%%%%%%%%%%%%%%%%%%%%%%%%%%%%%%%%%%%%%%%%%%%%%%%%%%%%%%%%%%%%%%%%%%%%%%%%%%%%%%%%%%%%%%%%%%
%%%%%%%%%%%%%%%%%%%%%%%%%%%%%%%%%%%%%%%%%%%%%%%%%%%%%%%%%%%%%%%%%%%%%%%%%%%%%%%%%%%%%%%%%%%%
%%%%%%%%%%%%%%%%%%%%%%%%%%%%%%%%%%%%%%%%%%%%%%%%%%%%%%%%%%%%%%%%%%%%%%%%%%%%%%%%%%%%%%%%%%%%
%%%%%%%%%%%%%%%%%%%%%%%%%%%%%%%%%%%%%%%%%%%%%%%%%%%%%%%%%%%%%%%%%%%%%%%%%%%%%%%%%%%%%%%%%%%%
%%%%%%%%%%%%%%%%%%%%%%%%%%%%%%%%%%%%%%%%%%%%%%%%%%%%%%%%%%%%%%%%%%%%%%%%%%%%%%%%%%%%%%%%%%%%
%%%%%%%%%%%%%%%%%%%%%%%%%%%%%%%%%%%%%%%%%%%%%%%%%%%%%%%%%%%%%%%%%%%%%%%%%%%%%%%%%%%%%%%%%%%%
%%%%%%%%%%%%%%%%%%%%%%%%%%%%%%%%%%%%%%%%%%%%%%%%%%%%%%%%%%%%%%%%%%%%%%%%%%%%%%%%%%%%%%%%%%%%
%%%%%%%%%%%%%%%%%%%%%%%%%%%%%%%%%%%%%%%%%%%%%%%%%%%%%%%%%%%%%%%%%%%%%%%%%%%%%%%%%%%%%%%%%%%%
%%%%%%%%%%%%%%%%%%%%%%%%%%%%%%%%%%%%%%%%%%%%%%%%%%%%%%%%%%%%%%%%%%%%%%%%%%%%%%%%%%%%%%%%%%%%
%%%%%%%%%%%%%%%%%%%%%%%%%%%%%%%%%%%%%%%%%%%%%%%%%%%%%%%%%%%%%%%%%%%%%%%%%%%%%%%%%%%%%%%%%%%%

%%%%%%%%%%%%%%%%%%%%%%%%%%%%%%%%%%%%%%%%%%%%%%%%%%%%%%%%%%%%%%%%%%%%%%%%%%%%%%%%%%%%%%%%%%%%

\section{One-well dynamics of the small jump component\label{s:onewell}}

\subsection{Exponential estimate for the small-jump component \label{s:exp}}

For $\rho>0$ and $0<\e\leq 1$ consider 
the decomposition $L=\xi^\e+\eta^\e$, where the L\'evy processes $\xi^\e$ and $\eta^\e$ have 
generating triplets $(d,\nu^\e_\xi,\mu)$ and $(0,\nu^\e_\eta,0)$ respectively with
\be
\nu^\e_\xi(\cdot)=
\nu\left(\cdot\cap \left[-\frac{1}{\e^{\rho}},\frac{1}{\e^{\rho}}\right]\backslash\{0\}\right),\qquad
\nu^\e_\eta(\cdot)=
\nu\left(\cdot\cap \mathbb{R}\backslash\left[-\frac{1}{\e^{\rho}},\frac{1}{\e^{\rho}}\right]\right).
\ee
The absolute value of jumps of the process $\e\xi^\e$ does not exceed $\e^{1-\rho}$. 

Thus the process $\xi^\e$ has a L\'evy measure with compact support, and 
the L\'evy measure $\nu_\eta^\e(\cdot)$ of $\eta^\e$ is finite. Denote
\be
\label{eq:chi}
\beta_\e=\nu_\eta^\e(\mathbb R)=
\int_{\mathbb{R}\backslash [-\frac{1}{\e^{\rho}},\frac{1}{\e^{\rho}}]}
\nu(dy)=H\left(\frac{1}{\e^{\rho}}\right).
\ee

Then, $\eta^\e$ is a compound Poisson process with intensity
$\beta_\e$, and jumps distributed according to the law
$\beta_\e^{-1}\nu_\eta^\e(\cdot)$.

Denote $\tau^\e_k$, $k\geq 0$, the jump times of $\eta^\e$ with $\tau^\e_0=0$.
Let $T^\e_k=\tau^\e_k-\tau^\e_{k-1}$ denote successive inter-jump periods, and
$W^\e_k=\eta^\e_{\tau^\e_k}-\eta^\e_{\tau^\e_k-}$
the jump heights of $\eta^\e$.
Then, the three processes $(T^\e_k)_{k\ge 1},$ $(W^\e_k)_{k\ge 1}$, and $(\xi^\e)_{t\geq 0}$ are independent. Moreover,
\begin{align}
&\P(T^\e_k\geq u)= \int_u^\infty \beta_\e e^{-\beta_\e s}\,ds=e^{-\beta_\e u},
\quad u\geq 0,\quad \text{ and } \quad
\E T^\e_k=\frac{1}{\beta_\e},\\
\label{eq:w}
&\P(W^\e_k< w)=\frac{\nu_\eta^\e(-\infty,w)}{\nu_\eta^\e(\mathbb{R})}
= \frac{1}{\beta_\e}\int_{(-\infty,w)}
 \I{|y|>\frac{1}{\e^{\rho}}} \nu(dy), \quad
w\in\mathbb{R}.
\end{align}
Between the arrival times of $\eta^\e$ the process $X^\e$ is driven by $\e\xi^\e$.
The next Lemma shows that on long time intervals $\e\xi^\e$ does not essentially deviate from zero.  
Hence the dynamics of the process $X^\e$ on the intervals 
between arrival times of the process $\eta^\e$ can be seen as a small random perturbation of the 
underlying deterministic trajectory.

\begin{lemma}
\label{l:1}
For any $\rho\in (0,1)$, any $\g\in (0,1-\rho)$ and 
$\theta\in (0,1-\rho-\g)$ there is  $p_0>0$ and $\e_0>0$ such that  the inequality
\be
\P(\sup_{t\in [0,1/\e^{\theta}]}|\e\xi_t^\e|\geq \e^\g)\leq \exp{\left( -1/\e^p \right) }
\ee
holds for all $0<\e\leq \e_0$ and $0<p\leq p_0$.
\end{lemma}

\begin{proof}
Let $\rho$, $\g$ and $\theta$ be as in the statement of Lemma.
Since
\ba
\label{eq:in}
\P(\sup_{t\in [0,1/\e^{\theta}]}|\e\xi^\e_t|\geq \e^\g)\leq
\P(\sup_{t\in [0,1/\e^{\theta}]}(\e\xi^\e_t)\geq \e^\g)
+\P(\inf_{t\in [0,1/\e^{\theta}]}(\e\xi^\e_t)\leq -\e^\g).
\ea
we have to estimate two summands. Let us consider the first.

The L\'evy measure of $\e\xi^\e$ has compact support, hence the process $\e\xi^\e$ has exponential
moments. Moreover, $\e\xi^\e_t-\E(\e\xi^\e_t)$ is a zero-mean martingale, so that
\ba
\E(\e\xi^\e_t)&=\e \mu t +\e t\int_{-1/\e^\rho}^{-1}y\nu(dy)+\e t\int_1^{1/\e^\rho}y\nu(dy),\\
|\E(\e\xi^\e_t)|&\leq \e t\left[  \mu +\int_{|y|\geq 1}y\nu(dy)\right]=\e tm.
\ea
Then Kolmogorov's inequality for exponential functions of martingales yields
\ba
\P\left(\sup_{t\in [0,1/\e^{\theta}]}(\e\xi^\e_t )\geq \e^\g\right)
&\leq \P\left(\sup_{t\in [0,1/\e^{\theta}]}(\e\xi^\e_t - \E(\e\xi^\e_t)   )
\geq \e^\g-\e^{1-\theta}m\right)  \\
&=\P\left(\sup_{t\in [0,1/\e^{\theta}]}e^{u(\e\xi^\e_t - \E(\e\xi^\e_t))}
\geq e^{u(\e^\g-\e^{1-\theta}m)}\right) \\
&\leq e^{-u(\e^\g-\e^{1-\theta}m)} 
\sup_{t\in [0,1/\e^{\theta}]}\E e^{u(\e\xi^\e_t - \E(\e\xi^\e_t))}\\
&\leq e^{-u(\e^\g-2\e^{1-\theta}m)} 
\sup_{t\in [0,1/\e^{\theta}]}\E e^{u\e\xi^\e_t}.
\ea
where the latter exponent can be derived from the L\'evy--Hinchin representation,
\ba
\E \exp{(u\e\xi^\e_t)}
=\exp\left\lbrace dt\frac{\e^2 u^2}{2}+\mu t\e u
+t\int_{0<|y|\leq 1/\e^\rho} (e^{u\e y}-1-u\e  y \I{|y|\leq 1}) \nu(dy)\right\rbrace.
\ea
Denote
\be
\phi(u,\e,t )=\ln{\E \exp{(u\e\xi^\e_t)}}+2u m\e^{1-\theta}-\e^\g u
\ee
and let $u=u(\e)=1/\e^c$ for $c=(1-\rho+\g)/2$.
We show that $\sup_{t\in[0,1/\e^\theta]}\phi(u(\e),\e,t)\to-\infty$ as a power of $\e$.  Indeed, since 
$0<c<1-\rho$, a straightforward calculation yields
\ba
\sup_{t\in[0,1/\e^\theta]}&\phi(u(\e),\e,t)\leq  d\frac{\e^{2-\theta-2c}}{2} +|\mu|\e^{1-\theta-c}
+\frac{1}{\e^\theta}\int_{|y|\leq 1} (e^{\e^{1-c}y}-1-\e^{1-c}y)\nu(dy)\\
&+\frac{1}{\e^\theta}  \int_1^{1/\e^\rho}(e^{\e^{1-c}y}-1 )\nu(dy)
+\frac{1}{\e^\theta}\int_{-1/\e^\rho}^{-1}(e^{\e^{1-c}y}-1 )\nu(dy)
+2 m\e^{1-\theta-c}-\e^{\g-c}\\
&\leq \e^{2-\theta-2c} \left(\frac{d}{2}+\int_{|y|\leq 1}y^2\nu(dy)\right)
+ 2\e^{1-c-\rho-\theta}\int_1^\infty \nu(dy) +(2 m+|\mu|)\e^{1-\theta-c} -\e^{\g-c}.
\ea 
Then since
$2-\theta-2c, 1-c-\rho-\theta, 1-\theta-c>\g-c$ we can take
%
% We can always choose some positive $c$ that the following inequalities hold
% \be
% \begin{cases}
% 2-2c-\theta>\g-c,\\
% 1-c-\rho-\theta>\g-c,\\
% 1-\theta-c>\g-c,\\
% c-\g>0,\\ 
% 0<c<1-\rho,\\
% 0<\rho<1,\\
% 0<\g<1-\rho,\\
% 0<\theta<1-\rho-\g.
% \end{cases}
% \quad
% \Leftrightarrow\quad
% \begin{cases}
% \g<c<2-\theta-\g,\\
% c<1-\rho,\\
% 0<\rho<1,\\
% 0<\g<1-\rho,\\
% 0<\theta<1-\rho-\g,
% \end{cases}
% \ee
% for example, we can take $c=1-\rho/2$.  
$p_0=(c-\g)/2$ to obtain
\be
\sup_{t\in[0,1/\e^\theta]}\phi(u(\e),\e,t)\leq - \frac{1}{\e^p},\quad \e\downarrow 0,
\ee
for all $0<p\leq p_0$.

The inequality for $\inf$ is proved analogously.
\end{proof}

\subsection{Dynamics on compact interval, $a>-\infty$}

Our goal is to study the one-well dynamics of the small-jump process $x^\e$ and its unperturbed conterpart $x^0$, 
\ba
\label{eq:eqs}
x_t^\e(x)&=x-\int_0^t U'(x_{s-}^\e(x))\, ds+\e\xi_t^\e,\\
x^0_t(x)&=x-\int_0^t U'(x^0_s(x))\, ds,\quad t\geq 0.
\ea
For definiteness we assume that the well's minimum is located at the origin and thus
the corresponding domain of attraction for $x^0$ is $(a,b)$, $-\infty<a<0<b<+\infty$, if the well is inner, 
and $(-\infty,b)$ if
it is peripheral. In the first case we also assume that $a$ and $b$ 
are non-degenerate local maxima of $U$. In the 
second case,  $b$ is a non-degenerate local maximum and $U'(x)$ increases to 
infinity faster than linearly as $x\to-\infty$.
Denote the critical point curvatures as
$U''(0)=M_0>0$, $U''(b)=M_b>0$ and $U''(a)=M_a<0$ (when defined).

For $\g>0$ and $t\geq 0$ we introduce an event
\be
\mathcal{E}_{t}=\{\omega:\sup_{s\in[0,t]}|\e\xi_s^\e|\leq \e^{4\gamma}\}.
\ee

We prove the following estimates.

\begin{prop}
\label{p:3}
For any $\gamma>0$, any $c>0$ there is $\e_0>0$ such that for
$0<\e\leq \e_0$ the inequality
\be
\sup_{s\in [0,t]}|x_s^\e(x)-x^0_s(x)|\leq c\e^{2\g}
\ee
holds a.s.\ on the event 
$\mathcal{E}_{t}$
uniformly for $t\geq 0$ and $x\in [a+\e^\g, b-\e^\g]$.

\end{prop}

Consider the representation of the process
$x^\e$ in powers of $\e$
\be
x_t^\e(x)=x^0_t(x)+\e Z_t^\e(x)+R_t^\e(x),\quad t\geq 0,
\ee
where $Z^\e$ is the first approximation of $x^\e$ satisfying the
stochastic differential equation
\be
\label{eq:Z}
Z_t^\e(x)=-\int_0^t U''(x^0_s(x)) Z_{s-}^\e(x)\, ds+\xi_t^\e
\ee
and the remainder $R^\e(x)$ is the absolutly continuous function starting at $0$
and satisfying the integral equation
\be
\label{eq:R}
R^\e_t(x)=\int_0^t
\left[-U'(x^0_s(x)+\e Z_{s-}^\e(x)+R^\e_s(x))+U'(x^0_s(x))+U''(x^0_s(x))\e Z_{s-}^\e(x)\right]\,ds.
\ee
We shall prove two Lemmas about the small noise dynamics of these processes.

\begin{lemma}
\label{l:z}
There is a universal constant $C_Z>0$ such that for any $\gamma>0$
there is $\e_0>0$ such that for $0<\e\leq \e_0$ the inequality 
\be
\sup_{s\in[0,t]}|\e Z^\e_s(x)| \leq C_Z \e^{3\g}
\ee
holds a.s.\ on the event $\mathcal{E}_t$
uniformly for $t\geq 0$ and $x\in [a+\e^\g,b-\e^\g]$.
\end{lemma}

\begin{lemma}
\label{l:r}
There is a universal 
constant $C_R>0$ such that for any $\g>0$ 
there is $\e_0>0$ such that for 
$0<\e\leq \e_0$ the inequality
\be
\sup_{s\in[0,t]}|R^\e_s(x)|\leq  C_R \e^{3\gamma}
\ee
holds a.s.\ on the event $\mathcal{E}_t$
uniformly for $t\geq 0$ and $x\in [a+\e^\g,b-\e^\g]$.
\end{lemma}

\noindent
\textbf{The Proof of Proposition~\ref{p:3}} follows easily from the previous Lemmas. 
 \hfill $\blacksquare$

\medskip

The proof of Lemmas \ref{l:z} and \ref{l:r} is performed in the sequel. 
We consider in detail only the neighbourhood of the critical point $a$. The behaviour of $x^\e$ in the 
neighbourhood of $b$ is obviously similar.
The following geometric properties of
the potential $U$ will be extensively used:
\begin{enumerate}
\item The deterministic trajectories $x^0_t(x)$, $x\in [a+\e^\g,b-\e^\g]$ converge to $0$
as $t\to\infty$
due to the property $x U'(x)>0$ for $x\neq a,b,0$.
\item The curvature of the potential at $x=a,b$ is negative. In a small
neighbourhood of $a$ we have
$U(x)=U(a)- M_a\frac{(x-a)^2}{2}+o((x-a)^2)$. Consequently $x^0$ behaves
there like
$a+e^{M_a t}$, and the dynamics of $x^\e$ reminds of the dynamics of an inverted
process of Ornstein-Uhlenbeck type. 

\item The curvature of the potential at $x=0$ is positive. In small
neighbourhoods of $0$ we have
$U(x)=U(0)+ M_0\frac{x^2}{2}+o(x^2)$. Consequently $x^0$ decays
there like
$e^{-M_0 t}$, and the dynamics of $x^\e$ reminds of the dynamics of a
process of Ornstein-Uhlenbeck type.
\end{enumerate}

From now on, let $\g>0$ be fixed.
Using assumptions on $U$, for technical reasons we fix some small $\delta$, $0<\delta<\min\{|a|,b\}$, 
and consider $\delta$-neighbourhoods of the 
critical points $a$, $0$ and $b$  with the following properties:
\begin{itemize}
\item 
there are some $0<m^a_1\leq M_a\leq m^a_2$, $\frac{m^a_2}{m^a_1}< \frac{3}{2}$, such that if $a\leq x\leq a+\delta$ then 
$m^a_1 (x-a)\leq -U'(x)\leq m^a_2 (x-a)$;
\item 
$-U'(\cdot)$ is  monotone increasing in $x\in[a,a+\delta]$.
\item Similar estimates hold in $\delta$-neighbourhood of $b$. 
\item
There are 
some $0<m^0_1<m^0_2$ such that the inequality
$m^0_1 < \inf_{|x|<\delta} U''(x)\leq \sup_{|x|<\delta} U''(x)< m^0_2$ holds. 
\end{itemize}

For $\e$ such that $0<\e^\g<\delta$ and for $x\in[a+\e^\g,b-\e^\g]$, denote 
\ba
t_\e(x)&=
\begin{cases}
\text{ the first time } x^0_t(x) \text{ reaches the level }  a+\delta 
\text{ if } x\in[a+\e^\g,a+\delta],\\
\text{ the first time } x^0_t(x) \text{ reaches the level }  b-\delta 
\text{ if } x\in[b-\delta,b-\e^\g],\\
0, \text{ if } x\in[a+\delta, b-\delta],
\end{cases}\\
&=
\begin{cases}
\int_x^{a+\delta}\frac{dy}{-U'(y)},
&\text{ if } x\in[a+\e^\g,a+\delta],\\
\int_{b-\delta}^x\frac{dy}{U'(y)},
&\text{ if } x\in[b-\delta,b-\e^\g],\\
0, &\text{ if } x\in[a+\delta, b-\delta].
\end{cases}
\ea
Also define the time period
\be
\hat{T}=\max\{\int_{a+\delta}^{-\delta}\frac{dy}{-U'(y)}, 
\int_\delta^{b-\delta}\frac{dy}{U'(y)}\}.
\ee
$\hat{T}$ has the property that for all $x\in [a+\delta,b-\delta]$ and $t\geq \hat{T}$,
$|x^0_t(x)|\leq \delta$, i.e.\ after $\hat{T}$ the trajectory of $x^0(x)$ is within a
$\delta$-neighborhood of the stable point $0$.

\subsubsection{Estimates on $Z^\e$}

\textbf{Proof of Lemma~\ref{l:z} }

The solution to equation~\eqref{eq:Z} is explicitly given by
\be
Z_t^\e(x)=\int_0^t e^{-\int_s^t U''(x^0_u(x))\, du}\,d\xi_s^\e.
\ee
Integration by parts results in the following representation for $Z^\e$:
\be
\label{eq:z3}
Z_t^\e(x)=\xi_t^\e-
\int_0^t \xi_{s-}^\e U''(x^0_s(x))e^{-\int_s^t U''(x^0_u(x))\, du} \,ds.
\ee
For $x=0$, $x^0_t(x)=0$ for all $t\geq 0$, and $Z^\e(0)$ is a process of the Ornstein-Uhlenbeck
type starting at zero and given by the equation
\be
Z_t^\e(0)=\xi_t^\e-
M_0\int_0^t \xi_{s-}^\e e^{-M_0(t-s)} \,ds,
\ee
and hence for any $t\geq 0$
\be
\sup_{s\in[0,t]}|Z_s^\e(0)|\leq 2\sup_{s\in[0,t]}|\xi_s^\e|.
\ee
Further, it follows from \eqref{eq:z3} that for $t\geq 0$ and $x\in[a+\e^\g,b-\e^\g]$
\ba
\label{eq:z1}
\sup_{s\in[0,t]}|Z_s^\e(x)|
\leq 
%\left( 1+\sup_{s\in[0,t]}\int_0^s |U''(x^0_u(x))|e^{-\int_u^s U''(x^0_v(x))\, dv} \,du\right)
%\sup_{s\in[0,t]}|\xi_s^\e|\\
%&=
\left( 1+\int_0^t |U''(x^0_s(x))|e^{-\int_s^t U''(x^0_u(x))\, du} \,ds\right)
\sup_{s\in[0,t]}|\xi_s^\e|. 
\ea

In order to prove Lemma~\ref{l:z},
we distinguish three cases:  $x\in [a+\delta,b-\delta]$, $x\in [a+\e^\g,a+\delta]$ and 
$x\in [b-\delta,b-\e^\g]$. 

\noindent
1.\ Let $x\in [a+\delta,b-\delta]$. Then we show that for any $t\geq 0$ and for some positive $C_1$ 
\be
\sup_{s\in[0,t]}|Z^\e_s(x)|\leq C_1 \sup_{s\in[0,t]}|\xi^\e_s|. 
\ee
Let
\be
C_2=\max_{x\in [a+\delta,b-\delta]}
\int_0^{\hat{T}} |U''(x^0_s(x))|e^{-\int_s^{\hat{T}} U''(x^0_u(x))\, du} \,ds.
\ee
Consider an arbitrary $t\geq \hat{T}$. Then
\ba
\label{eq:2summand}
\int_0^t |U''(x^0_s(x))&|e^{-\int_s^t U''(x^0_u(x))\, du} \,ds\\
&=\int_0^{\hat{T}} |U''(x^0_s(x))|e^{-\int_s^t U''(x^0_u(x))\, du} \,ds
+\int_{\hat{T}}^t |U''(x^0_s(x))|e^{-\int_s^t U''(x^0_u(x))\, du} \,ds.
\ea
Let us estimate the first summand in \eqref{eq:2summand}. Since for all $x\in [a+\delta,b-\delta]$
and $t\geq \hat{T}$, $m^0_1 < U''(x^0_t(x))< m^0_2$, we have
\ba
\int_0^{\hat{T}} |U''(x^0_s(x))|e^{-\int_s^t U''(x^0_u(x))\, du} \,ds
&=e^{-\int_{\hat{T}}^t U''(x^0_u(x))\, du}
\int_0^{\hat{T}} |U''(x^0_s(x))|e^{-\int_s^{\hat{T}} U''(x^0_u(x))\, du} \,ds\\
&\leq e^{-m^0_1(t-\hat{T})}C_2\leq C_2.
\ea
The second summand in \eqref{eq:2summand} is estimated analogously:
\ba
\int_{\hat{T}}^t |U''(x^0_s(x))|e^{-\int_s^t U''(x^0_u(x))\, du} \,ds\leq
m^0_2\int_{\hat{T}}^t e^{-m^0_1(t-s)}\,ds\leq \frac{m^0_2}{m^0_1}.
\ea
Taking $C_1=\max\{2, C_2+\frac{m^0_2}{m^0_1}\}$ completes the proof.

\noindent
2.\ Let $x\in [a+\e^\g,a+\delta]$. Then we show that
\be
\sup_{s\in[0,t]}|Z^\e_s(x)|\leq \frac{C_3}{\e^\gamma} \sup_{s\in[0,t]}|\xi^\e_s|. 
\ee

Indeed, for $x\in[a+\e^\g, a+\delta]$ and $t\leq t_\e(x)$ we have,
\ba
\label{eq:z2}
1&+\int_0^t |U''(x^0_s(x))|e^{-\int_s^t U''(x^0_u(x))\, du} \,ds
=1-\int_0^t U''(x^0_s(x))e^{-\int_s^t U''(x^0_u(x))\, du} \,ds\\
&=1-\int_0^t U''(x^0_s(x))e^{\int_{x^0_s(x)}^{x^0_t(x)}\frac{ U''(v)}{ U'(v)}\, dv} \,ds
\quad (v=x^0_u(x), dv=-U'(v)du)\\
&=1-\int_0^t U''(x^0_s(x))e^{\ln{U'(x^0_t(x))/}U'(x^0_s(x))} \,ds\\
&=1-U'(x^0_t(x))\int_0^t \frac{U''(x^0_s(x))}{U'(x^0_s(x))} \,ds\\
&=1+U'(x^0_t(x))\int_x^{x^0_t} \frac{U''(v)}{U'(v)^2} \,dv
\quad (v=x^0_s(x), dv=-U'(v)ds)\\
&=1-U'(x^0_t(x))\left(\frac{1}{U'(x^0_t(x))}-\frac{1}{U'(x)} \right) 
=\frac{U'(x^0_t(x))}{U'(x)}.
\ea

For any $t\geq 0$ we use \eqref{eq:z1} and \eqref{eq:z2} to obtain  
\ba
&1+\int_0^t |U''(x^0_s(x))|e^{-\int_s^t U''(x^0_u(x))\, du} \,ds\\
&=1- e^{-\int_{t_\e(x)\wedge t}^t U''(x^0_u(x))\, du}   
\int_0^{t_\e(x)\wedge t} 
U''(x^0_s(x)) e^{-\int_s^{t_\e(x)\wedge e} U''(x^0_u(x))\, du} \,ds\\
&+\int_{t_\e(x)\wedge t}^t |U''(x^0_s(x))|e^{-\int_s^t U''(x^0_u(x))\, du} \,ds\\
&=1-\frac{U'(x^0_t(x))}{U'(x^0_{t_\e(x)\wedge t}(x))}
\left(1-\frac{U'(x^0_{t_\e(x)\wedge t}(x))}{U'(x)} \right) \\
&+ \int_0^{t-t_\e(x)\wedge t} 
|U''(x^0_s(a+\delta))|e^{-\int_s^{t-t_\e(x)\wedge t} U''(x^0_u(a+\delta))\, du} \,ds\\
&\leq 1-\frac{U'(x^0_t(x))}{U'(x^0_{t_\e(x)\wedge t}(x))}+\frac{U'(x^0_t(x))}{U'(x)}+C_1.
\ea
Note that for  $0< \e\leq \e_0$ for some $\e_0$ small enough and depending on $U$, $a$, $b$, $\g$
and $\delta$ 
\ba
\label{eq:j}
1-\frac{U'(x^0_t(x))}{U'(x^0_{t_\e(x)\wedge t}(x))}+\frac{U'(x^0_t(x))}{U'(x)}\leq
\begin{cases}
\frac{U'(x^0_t(x))}{U'(x)}\leq \frac{m_2^a|a+\delta|}{m_1^a(a-x)}\leq \frac{C_3}{\e^\g}, 
\quad 0\leq t\leq t_\e(x),\\
1+\max_{y\in[a,b]}|U'(y)|(\frac{1}{|U'(a+\delta)|}+\frac{1}{|U'(x)|})
\leq \frac{C_4}{\e^\g}, \quad  t\geq t_\e(x),\\
\end{cases}
\ea

\noindent
3.\ For $x\in [b-\delta,b-\e^\g]$ we obtain an estimate similar to \eqref{eq:j}.

Hence, for all $x\in[a+\e^\g,b-\e^\g]$, $t\geq 0$ and $0<\e\leq\e_0$,
\be
\sup_{s\in[0,t]}|\e Z^\e(x)|\leq \frac{C_Z}{\e^\gamma} \sup_{s\in[0,t]}|\e\xi^\e_s|\leq 
 C_Z \e^{3\gamma}
\ee
for some positive $C_Z$ on the event 
$\mathcal{E}_t$. \hfill $\blacksquare$

%%%%%%%%%%%%
%%%%%%%%%%%%
%%%%%%%%%%%%
%%%%%%%%%%%%
%%%%%%%%%%%%
%%%%%%%%%%%%

\subsubsection{Estimates on $R^\e$}

To estimate the remainder term $R^\e$ we need finer smoothness
properties of the potential $U$. However, the following Lemma shows that this restriction only has
to hold locally.

\begin{lemma}
\label{l:Rrough}
There exists $C>0$ and $\e_0>0$ such that for and $0<\e\leq\e_0$ the inequality 
\be
\sup_{s\in[0,t]}|R^\e_s(x)|\leq C
\ee
holds a.s.\ on the event $\mathcal{E}_t$ uniformly for $t\geq 0$ and $x\in [a+\e^\g,b-\e^\g]$.

\end{lemma}
\begin{proof}
By Assumption \textbf{U2} we know that for any $t\geq 0$, $x\in [a,b]$ we have 
$x^0_t(x)\in [a,b]$.
Moreover, for $0<\e\leq \e_Z$ and $x\in[a+\e^\g,b-\e^\g]$ we have
$\sup_{s\in [0,t]}|\e Z_s^\e(x)| < C_Z\e^{3\gamma}$ on $\mathcal{E}_t$ due to Lemma \ref{l:z}.
Recall that $U'$ increases at least linearly at infinity (see Assumption \textbf{U3}). This guarantees the 
existence of a constant
$C > 0$ such that for any $x\in [a,b]$, $|z|\leq 1$  we have
\be
-U'(x+z+C)+ U'(x)+ U''(x)z<0.
\ee
Hence for any $0\leq s\leq t$, $x\in[a+\e^\g,b-\e^\g]$ the inequality
\be
-U'(x^0_s(x)+ \e Z_{s-}^\e(x)+C)+ U'(x^0_s(x))+ U''(x^0_s(x))\e Z_{s-}^\e(x)<0
\ee
holds on the event $\mathcal{E}_t$ for 
$0<\e\leq \min\{C_Z^{-1/3\g},\e_Z\}$. Now assume there is some $x\in[a+\e^\g,b-\e^\g]$,
and some (smallest) $\tau\in [0,t]$
such that $R^\e_\tau(x)=C$. Observe that the rest term $R^\e$ satisfies the integral equation
\be
R^\e_s(x)=\int_{0}^s f(R^\e_u(x), x^0_u(x), \e Z^\e_{u-}(x))\,du
\ee
with the smooth integrand
\ben
f(R, x^0, \e Z)=-U'(x^0+\e Z+R)+U'(x^0)+U''(x^0)(\e Z).
\een
This implicitly says that $R^\e$ is an absolutely 
continuous function of time. By definition of $\tau$,
we have
\be 
0 \le D^+R_s^{\e}(x)\vert_{s=\tau} = 
-U'(x^0_{\tau}(x)+ \e Z_{\tau-}^\e(x)+C)+ U'(x^0_{\tau}(x))+ U''(x^0_{\tau}(x))\e Z_{\tau-}^\e(x)<0,
\ee
a contradiction, with $D^+$ denoting the right Dini derivative. A similar reasoning applies under the assumption
$R^\e_\tau(x)= - C$. This completes the proof.

\end{proof}

Lemma \ref{l:Rrough} has a very convenient consequence. It states precisely that the solution 
process
$x^\e_s(x)$, $s\in[0,t]$, with initial state confined to $[a+\e^\g,b-\e^\g]$, 
stays bounded by a deterministic constant $K$ on 
the set $\mathcal{E}_t$, $t\geq 0$. 
Therefore, in the small noise limit,
only local properties of $U$ are relevant to our analysis.

\begin{lemma}
\label{l:r1}
There exists $C_1>0$ such that for any $\gamma>0$ there is 
$\e_0>0$ such that for $0\leq \e\leq \e_0$,  
\be
\label{eq:r1}
\sup_{s\in [0,t\wedge t_\e(x)]}|R^\e_s(x)|
\leq  C_1 \e^{3\g}
\ee
on the event $\mathcal{E}_t$ uniformly for $x\in [a+\e^\g,b-\e^\g]$ and $t\geq 0$.
\end{lemma}

\begin{proof}
\noindent
1.\
For $x\in[a+\delta, b-\delta]$ the time $t_\e(x)=0$ and the estimate \eqref{eq:r1} is trivial. 
Thus it is only necessary to consider $x$ from the neighbourhoods of the boundary points $a$ and $b$. 
For definiteness, we consider the case $x\in[a+\e^\g, a+\delta]$.
Let also Lemmas \ref{l:z} and \ref{l:Rrough} hold for $0<\e\leq\e_1$ with constants $C_Z$ and $C$.

\noindent
2.\ The rest term $R^\e$ satisfies the integral equation
\be
\label{eq:r2}
R^\e_t(x)=\int_0^t f(R^\e_s(x), x^0_s(x), \e Z^\e_{s-}(x))\,ds
\ee
with
\be
f(R, x^0, \e Z)=-U'(x^0+\e Z+R)+U'(x^0)+U''(x^0)(\e Z).
\ee
Moreover, $R^\e$ is an absolutely continuous function of time.
Let the constant $K$ from Assumption \textbf{U1} be bigger than $C$.
We write the Taylor expansion for the integrand $f$ with some $|\theta|\leq K$:
\ba
f(R, x^0, \e Z)&=-U'(x^0+\e Z+R)+U'(x^0)+U''(x^0)(\e Z)\\
&=-U'(x^0)-U''(x^0)(R+\e Z)-\frac{U^{(3)}(\theta)}{2}(R+\e Z)^2+U'(x^0)+ U''(x^0)(\e Z) \\
&=-U''(x^0)R-\frac{U^{(3)}(\theta)}{2}(R+\e Z)^2
\ea
Since $U\in {\cal C}^3$, $|U^{(3)}|$ is bounded, say by $L$, on $[-K,K]$. Using the
inequality $(R+\e Z)^2\leq 2(R^2+\e^2 Z^2)$ we
obtain that for $t\geq 0$, 
\ba
f(R^\e_t(x),x^0_t(x),\e Z_{t-}^\e(x))&\leq -U''(x^0_t(x))R^\e_t +L (R^\e_t)^2+L(\e Z^\e_{t-}(x))^2 
< -U''(x^0_t(x))R^\e_t +L (R^\e_t)^2 + A^2\e^{6\gamma},\\
f(R^\e_t(x),x^0_t(x),\e Z_{t-}^\e(x))&\geq -U''(x^0_t(x))R^\e_t -L (R^\e_t)^2-L(\e Z^\e_{t-}(x))^2 
> -U''(x^0_t(x))R^\e_t -L (R^\e_t)^2 - A^2\e^{6\gamma}
\ea
on the event 
$\mathcal{E}_t$,  with $A^2=2C_Z^2L$.

\noindent
3.\ Let us prove the upper bound in \eqref{eq:r1}.
Together with \eqref{eq:r2} consider the Riccati equation
\be
p^\e_t(x)=\int_0^t (m_2^a p^\e_s +L (p_s^\e)^2 + A^2\e^{6\g})\,ds,
\quad 0\leq t\leq t_\e(x).
\ee
Under the conditions of the lemma, it is enough to prove two statements:
\begin{description}
\item[a)] $R^\e_t(x)\leq p^\e_t$ for $0\leq t\leq t_\e(x)$.
\item[b)] $p^\e_t\leq  C_1\e^{3\gamma}$ for $0\leq t\leq t_\e(x)$.
\end{description}
We have the closed form formula for $p_t$: 
\ba
p^\e_t&=  A^2\e^{6\g} \frac{ e^{t \lambda^\e  }-e^{-t \lambda^\e  } }
{(m^a_2+ \lambda^\e         )e^{-t \lambda^\e  }
-(m^a_2- \lambda^\e   )e^{t   \lambda^\e   }},\\
\lambda^\e&=\sqrt{(m^a_2)^2-4LA^2\e^{6\g}}.
\ea
It is easy to see that $p^\e_t$ is a non-negative monotonically increasing function starting at $0$.
However $p^\e_t$ has a singularity at 
\be
t^\ast(\e)=\frac{1}{2   \lambda^\e         } \ln\left( 
\frac{m^a_2+ \lambda^\e        }
{m^a_2-    \lambda^\e           }\right) \geq
\frac{3\gamma|\ln\e|}{m^a_2}(1+\mathcal{O}(|\ln\e|^{-1})),
\ee
where the latter inequality holds for $\e\downarrow 0$.
Note that
\be
\label{eq:tt}
t_\e(x)=\int_{x}^{a+\delta}\frac{dy}{|U'(y)|}
\leq \frac{1}{m^a_1}\int_{a+\e^\g}^{a+\delta}\frac{dy}{|a-y|}
=\frac{1}{m^a_1}\ln{\left( \frac{\delta}{\e^\g}\right) }=t_\e.
\ee
In the limit of small $\e$, $t_\e$ can be calculated as
\be
t_\e=\frac{\g|\ln\e|}{m^a_1}(1+\mathcal{O}(|\ln\e|^{-1})),
\ee
Hence $t_\e(x)\leq t_\e< t^\ast(\e)$ for $0<\e\leq \e_2$, $\e_2$ being sufficiently small, 
and $p^\e_t$ is well defined on the time interval under 
consideration.

To show \textbf{a)} we note that at the starting point $t=0$,
\be
D^+R^\e_t(x)\Big|_{t=0}=\lim_{h\downarrow 0}\frac{R^\e_h(x)-0}{h}
=0<A^2\e^{6\g}=\dot p^\e_t(x)\Big|_{t=0},
\ee
consequently it follows from the continuity of $R^\e$ and $p^\e$ that $p^\e_t>R^\e_t$ for
at least positive and small $t$.
Assume there exists $\tau=\inf\{t>0: R^\e_\tau(x)=p^\e_\tau\}$ such that $\tau\leq t_\e(x)$.
At the point $\tau$ the left derivative of $R^\e(x)$ is necessarily not less than the
derivative of $p^\e$ which leads to the following contradiction:
\ba
D^-R^\e_t(x)\Big|_{t=\tau}&=\lim_{h\downarrow 0}\frac{R^\e_\tau(x)-R^\e_{\tau-h}(x)}{h}
=f(R^\e_\tau(x),x^0_\tau(x),Z^\e_{\tau-}(x))\\
&\geq \dot p^\e_t\Big|_{t=\tau}=m_2^a p^\e_\tau +L (p^\e_\tau)^2 + A^2\e^{6\g},\\
f(R^\e_\tau(x),x^0_\tau(x),Z^\e_{\tau-}(x))&=f(p^\e_\tau,x^0_\tau(x),Z^\e_{\tau-}(x))
<m_2^a p^\e_\tau +L (p^\e_\tau)^2 + A^2\e^{6\g}.
\ea
To prove \textbf{b)}, we use the inequality
\be
\sup_{t\in [0,t_\e(x)]}p^\e_t\leq p^\e_{t_\e(x)}\leq p^\e_{t_\e}.
\ee
and a formula \eqref{eq:tt} for $t_\e$.
Indeed, on $\mathcal{E}_t$, we have the following estimates
\ba
&e^{t_\e  \lambda^\e}-e^{-t_\e  \lambda^\e  }
\leq\left(\frac{\delta}{\e^\g}\right)^\frac{ \lambda^\e  }{m_1^a}
\leq\left(\frac{\delta}{\e^\gamma}\right)^\frac{m_2^a}{m_1^a}=c_1\e^{-\g\frac{m_2^a}{m_1^a}},\\
&(m_2^a+   \lambda^\e       )e^{-t_\e   \lambda^\e   }\geq 
m_2^a\left(\frac{\delta}{\e^\g}\right)^{-\frac{   \lambda^\e   }{m_1}}
\geq c_2\e^{\g\frac{m_2^a}{m_1^a}},\\
&(m_2^a-   \lambda^\e         )e^{t_\e   \lambda^\e   }
\leq m_2^a\left(1-\sqrt{1-\frac{4LA^2\e^{6\gamma}}{(m_2^a)^2}} \right) 
\left(\frac{\delta}{\e^\gamma}\right)^\frac{   \lambda^\e   }{m_1^a}\\
&\qquad\leq \frac{4LA^2\e^{6\gamma}}{m_2}
\left(\frac{\delta}{\e^\gamma}\right)^\frac{m_2^a}{m_1^a}
=c_3\e^{\gamma(6-\frac{m^a_2}{m^a_1})}.
\ea
Thus, since $\frac{m_2^a}{m_1^a}<\frac{3}{2}$ and for $\e\leq \e_0=\min\{\e_1,\e_2\}$ we can estimate
\be
p^\e_{t_\e}\leq A^2\e^{6\gamma} \frac{c_1\e^{-\gamma\frac{m_2^a}{m_1^a}}}
{c_2\e^{\gamma\frac{m_2^a}{m_1^a}}- c_3\e^{\gamma(6-\frac{m_2^a}{m_1^a})}}
=\e^{2\gamma(3-\frac{m^a_2}{m^a_1})} \frac{A^2c_1 }
{c_2-c_3\e^{2\gamma(3-\frac{m_2^a}{m_1^a})}}
\leq  C_1 \e^{3\gamma}.
\ee

The proof of the lower bound in \eqref{eq:r1} is analogous.
\end{proof}

%%%%%%%%%%%
%%%%%%%%%%%
%%%%%%%%%%%%%%%%%%%%%%
%%%%%%%%%%%

\begin{lemma}[Estimate away from critical points]
\label{l:r2}
There exists $C_2>0$ such that for any $\gamma>0$ there is 
$\e_0>0$ such that for $0\leq \e\leq \e_0$ and any  $t_\e(x)\leq t\leq t_\e(x)+\hat{T}$,
\be
\sup_{s\in [t_\e(x),t]}|R^\e_s(x)|\leq C_2 \e^{3\gamma},
\ee
on the event $\mathcal{E}_t$
uniformly for $x\in[a+\e^\g,b-\e^\g]$.
\end{lemma}

\begin{proof}
Using Lemma \ref{l:Rrough}, choose $K>0$ such that on the event $\mathcal{E}_t$
the processes $x^\e(x), \e Z^{\e}(x), R^\e(x)$ take their values in $[-K, K]$ 
as long as time runs in $[0,t]$. Let also previous Lemmas hold for $0<\e\leq \e_0$.

For $t_\e(x)\leq t\leq t_\e(x)+\hat{T}$, the rest term $R^\e$ 
satisfies the following integral equation:
\ba
R^\e_t(x)&=R^\e_{t_\e(x)}(x)+\int_{t_{\e}(x)}^t
\left[-U'(x^0_s(x)+\e Z_{s-}^\e(x)+R^\e_s(x))+U'(x^0_s(x))+U''(x^0_s(x))\e Z_{s-}^\e(x)\right]\, ds\\
&=R^\e_{t_\e(x)}(x)
-\int_{t_\e(x)}^t \left[U'(x^0_s(x)+\e Z_{s-}^\e(x)+R^\e_s(x))-U'(x^0_s(x)+\e Z_{s-}^\e(x)) \right]\, ds\\
&-\int_{t_\e(x)}^t \left[U'(x^0_s(x)+\e Z_{s-}^\e(x))-U'(x^0_s(x))-U''(x^0_s(x))\e Z_{s-}^\e(x)  \right]\, ds\\
&=R^\e_{t_\e(x)}(x)-\int_{t_\e(x)}^t U''(\theta^1_{s})R^\e_s(x)  \, ds
-\int_{t_\e(x)}^t \tfrac{1}{2}U^{(3)}(\theta^2_{s})(\e Z_{s-}^\e(x))^2  \, ds
\ea
with appropriate $\theta^1_s$, $\theta^2_s\in [-K,K]$.
Note that $R^\e_{t_\e(x)}(x)=0$  if $x\in[a+\delta,b-\delta]$.

Thus on $\mathcal{E}_t$, with the help of Lemma~\ref{l:r1}, we obtain
\ba
|R_t^\e(x)|&\leq |R^\e_{t_\e(x)}(x)|+\int_{t_\e(x)}^t L|R^\e_s(x)|  \, ds +
\tfrac{t-t_\e(x)}{2}LC^2_Z\e^{6\gamma}
\leq C_1\e^{3\gamma}+\int_{t_\e(x)}^{t_\e(x)+\hat{T}} L|R^\e_s(x)|  \, ds +
\tfrac{1}{2}\hat{T}LC^2_Z\e^{6\gamma}.
\ea
An application of Gronwall's lemma yields the final estimates for 
$t_\e(x)\leq t\leq t_\e(x)+\hat{T}$:
\ba
|R_t^\e(x)|&\leq \left( C_1\e^{3\gamma}+\tfrac{1}{2}\hat{T}LC^2_Z\e^{6\gamma}\right) 
e^{\hat{T}L}\leq C_2\e^{3\gamma}.
\ea
\end{proof}

%%%%%%
%%%%%%
%%%%%%
%%%%%%

\begin{lemma}[Estimate near the stable point]
\label{l:r3}
There exist a positive constant $C_3$ 
such that for any $\gamma>0$ there is 
$\e_0>0$ such that for $0\leq \e\leq \e_0$ and any $t\geq t_\e(x)+\hat{T}$,
\be
\sup_{s\in [t_\e(x)+\hat{T},t]}|R^\e_s(x)|\leq   
\e^{3\gamma},
\ee
on the event $\mathcal{E}_t$
uniformly for $x\in[a+\e^\g,b-\e^\g]$.
\end{lemma}

\begin{proof}
\noindent 1.\
Using Lemma \ref{l:Rrough}, choose $K>0$ such that on the event $\mathcal{E}_t$ the
processes $x^\e(x), \e Z^{\e}(x), R^\e(x)$ take their values in $[-K, K]$ as 
long as time runs in $[0,t]$. Let previous Lemmas hold for $0<\e\leq \e_0$.

For $t\geq t_\e(x)+\hat{T}$
the rest term $R^\e$ satisfies the integral equation
\be
R^\e_t(x)=R^\e_{t_\e(x)+\hat{T}}(x)
+\int_{t_\e(x)+\hat{T}}^t f(R^\e_s(x), x^0_s(x), \e Z^\e_{s-}(x))\,ds
\ee
with
\be
f(R, x^0, \e Z)=-U'(x^0+\e Z+R)+U'(x^0)+U''(x^0)(\e Z).
\ee
Note that for the time instants $t$ under consideration, the deterministic trajectory $x^0_t(x)$
is in the $\delta$-neighbourhood of the stable point $0$.
Repeating the argument of Lemma~\ref{l:r1} we obtain the following estimates: 
\ba
f(R^\e_t(x),&x^0_t(x),\e Z_{t-}^\e(x))\leq -U''(x^0_t(x))R^\e_t +L (R^\e_t)^2+L(\e Z^\e_{t-}(x))^2 \\
&< -U''(x^0_t(x))R^\e_t +L (R^\e_t)^2 + 2C_Z^2L^2\e^{3\gamma}<
-U''(x^0_t(x))R^\e_t +L (R^\e_t)^2 + D \e^{3\gamma},\\
f(R^\e_t(x),&x^0_t(x),\e Z_{t-}^\e(x))\geq -U''(x^0_t(x))R^\e_t -L (R^\e_t)^2-L(\e Z^\e_{t-}(x))^2\\ 
&> -U''(x^0_t(x))R^\e_t -L (R^\e_t)^2 - 2C_Z^2L^2\e^{3\gamma}
> -U''(x^0_t(x))R^\e_t -L (R^\e_t)^2 - D\e^{3\gamma}
\ea
on the event 
$\mathcal{E}_t$, 
with some $D>2C_Z^2L$ which will be specified later.

The main difference to Lemma~\ref{l:r1} consists in the sign of the $U''$ in the vicinity
of zero. Now the curvature is positive what guarantees the boundedness of $R^\e(x)$ on long
time intervals.

2.\ We establish the upper bound for $R^\e(x)$.
Consider a Riccati equation 
\be
\label{eq:ric}
p^\e_t=R^\e_{t_\e(x)+\hat{T}}(x)
+\int_{t_\e(x)+\hat{T}}^t (-m_2^0 p^\e_s +L (p^\e_s)^2 + D\e^{3\gamma})\, ds,
\quad t\geq t_\e(x)+\hat{T}.
\ee  

The comparison argument of Lemma~\ref{l:r1} shows that
\be
R^\e_t(x)\leq p^\e_t,\quad t\geq t_\e(x)+\hat{T}.
\ee
Now we study the Riccati equation \eqref{eq:ric} in detail. It is easy to see that it
has two positive stationary solutions at which the integrand of \eqref{eq:ric} vanishes:
\be
p^\pm=\frac{m_2^0}{2L}\left(1\pm\sqrt{1-\frac{4LD^2\e^{3\gamma}}{(m_2^0)^2}} \right). 
\ee
Applying the elementary inequality $\frac{x}{2}\leq 1-\sqrt{1-x}\leq x$, 
$x\in[0,1]$, to the smaller solution 
$p^-$ and for $\e\leq \e_0$ such that $4LD^2\e^{3\g}/(m_2^0)^1\leq 1$,  we find that
\be
\frac{D^2}{m_2^0}\e^{3\gamma}\leq p^- \leq \frac{2D^2}{m_2^0}\e^{3\g}
<\frac{2D^2}{m_1^0}\e^{3\gamma}.
\ee
 
This means that if $R_{t_\e(x)+\hat{T}}^\e(x)<\frac{D^2}{m_2^0}\e^{3\g}$, the solution 
$R^\e_t(x)$ does not exceed $\frac{2D^2}{m_1^0}\e^{3\g}$ on the time interval
$[t_\e(x)+\hat{T}, t]$ and the event $\mathcal{E}_t$.

We use Lemma~\ref{l:r2} to conclude that $R_{t_\e(x)+\hat{T}}^\e(x)<C_2\e^{3\g}$, 
and taking $D>\sqrt{C_2 m_2^0}$ and $C_3>\frac{2D^2}{m_1^0}$ finishes the proof.

3.\ The lower bound for $R^\e(x)$ is obtained analogously.

\end{proof}

\textbf{Proof of Lemma~\ref{l:r}} 
The claim of Lemma~\ref{l:r} follows from Lemmas~\ref{l:r1}, \ref{l:r2} and \ref{l:r3}
by taking $C_R=\max\{C_1,C_2,C_3\}$ and $\e_0$ the minimal value of $\e$ for which these Lemmas hold
simultaneously.
\hfill $\blacksquare$

%%%%%%%%%%%%%%%%%%%%%%%%%%%%%%%%%%%%%%%%%%%%%%%%%%%%%%%%%%%%%%%%%%%%%%%%%%%%%%%%%%%%%%%%%%%%%%%
%%%%%%%%%%%%%%%%%%%%%%%%%%%%%%%%%%%%%%%%%%%%%%%%%%%%%%%%%%%%%%%%%%%%%%%%%%%%%%%%%%%%%%%%%%%%%%%
%%%%%%%%%%%%%%%%%%%%%%%%%%%%%%%%%%%%%%%%%%%%%%%%%%%%%%%%%%%%%%%%%%%%%%%%%%%%%%%%%%%%%%%%%%%%%%%
%%%%%%%%%%%%%%%%%%%%%%%%%%%%%%%%%%%%%%%%%%%%%%%%%%%%%%%%%%%%%%%%%%%%%%%%%%%%%%%%%%%%%%%%%%%%%%%
%%%%%%%%%%%%%%%%%%%%%%%%%%%%%%%%%%%%%%%%%%%%%%%%%%%%%%%%%%%%%%%%%%%%%%%%%%%%%%%%%%%%%%%%%%%%%%%
%%%%%%%%%%%%%%%%%%%%%%%%%%%%%%%%%%%%%%%%%%%%%%%%%%%%%%%%%%%%%%%%%%%%%%%%%%%%%%%%%%%%%%%%%%%%%%%
%%%%%%%%%%%%%%%%%%%%%%%%%%%%%%%%%%%%%%%%%%%%%%%%%%%%%%%%%%%%%%%%%%%%%%%%%%%%%%%%%%%%%%%%%%%%%%%
%%%%%%%%%%%%%%%%%%%%%%%%%%%%%%%%%%%%%%%%%%%%%%%%%%%%%%%%%%%%%%%%%%%%%%%%%%%%%%%%%%%%%%%%%%%%%%%
%%%%%%%%%%%%%%%%%%%%%%%%%%%%%%%%%%%%%%%%%%%%%%%%%%%%%%%%%%%%%%%%%%%%%%%%%%%%%%%%%%%%%%%%%%%%%%%
%%%%%%%%%%%%%%%%%%%%%%%%%%%%%%%%%%%%%%%%%%%%%%%%%%%%%%%%%%%%%%%%%%%%%%%%%%%%%%%%%%%%%%%%%%%%%%%

\subsection{Final estimate for $|x^\e_t-x^0_t|$, $a>-\infty$}

In this section we use Lemma \ref{l:1} and Proposition \ref{p:3} to estimate the probability
that the small-jump process $x^\e_t(x)$ leaves the $\frac{1}{2}\e^{2\g}$-dependent tube of the 
deterministic trajectory $x^0_t(x)$.

\begin{prop}
\label{p:11}
Let $a>-\infty$. Let $\rho\in(0,1)$,  $x^\e(x)$ and $x^0(x)$ satisfy \eqref{eq:eqs}. 
Let $T(\e)$ be an exponentially distributed random variable with mean $1/\beta_\e$
and let $\xi^\e$ and let $T(\e)$ be independent. Then for any $\g\in (0,(1-\rho)/4)$
there exist $p_0>0$ and $\e_0>0$ such that the inequality
\ba
\sup_{x\in [a+\e^\g,b-\e^\g]}
\Px(\sup_{t\in [0,T(\e)]}|x^\e_t(x)-x^0_t(x)|\geq \frac{\e^{2\g}}{2})\leq \exp{\left( -1/\e^p \right) }
\ea
holds for all $0\leq p\leq p_0$ and $0<\e\leq \e_0$.
\end{prop}

\begin{proof} 
Let $\theta<\min{1-\rho-\g, r\rho}$. 
% Then according to Lemma \ref{l:1}
% $\P(\mathcal{E}^c_{1/\e^\theta})\leq \exp(-\e^{-p})$ for $0< p\leq p_1$ and $0<\e\leq \e_1$.
 Then
$\E T_1= \frac{1}{\beta_e}\gg \frac{1}{\e^\theta}$ as $\e\downarrow 0$. 
Consider the number 
\be
k_\e=\left[ \frac{\e^{\theta/2}}{\beta_\e}\right],
\ee
where $[x]$ denotes the integer part of $x$. Note that $k_\e\to\infty$ slower than some power of $1/\e$.

For any $x\in[a+\e^\g,b-\e^\g]$ we have
\ba
\Px&\left(\sup_{t\in[0,T_1]}|x^\e_t(x)-x^0_t(x)|\geq \frac{\e^{2\g}}{2}\right)=\int_0^\infty \beta_\e e^{-\beta_\e \tau}
\Px\left( \sup_{t\in[0,\tau]}|x^\e_t(x)-x^0_t(x)|\geq \frac{\e^{2\g}}{2}\right) \,d\tau\\
&=\left( \int_0^{k_\e/\e^\theta}+ \int_{k_\e/\e^\theta}^\infty\right) \beta_\e e^{-\beta_\e \tau}
\Px\left( \sup_{t\in[0,\tau]}|x^\e_t(x)-x^0_t(x)|\geq \frac{\e^{2\g}}{2}\right) \,d\tau. 
\ea
For $0<\e\leq \e_1$, $\e_1$ small enough, the second summand is estimated as
\ba
\label{eq:f}
\int_{k_\e/\e^\theta}^\infty \beta_\e e^{-\beta_\e \tau}
\Px\left( \sup_{t\in[0,\tau]}|x^\e_t(x)-x^0_t(x)|\geq \frac{\e^{2\g}}{2}\right) \,d\tau 
\leq
\exp{\left( -\frac{\beta_\e}{\e^\theta}k_\e \right) }
\leq \exp{\left( -\frac{1}{\e^{\theta/2}}\right) }.
\ea

For the first summand,
\ba
\int_0^{k_\e/\e^\theta}
&\beta_\e e^{-\beta_\e \tau}
\P\left( \sup_{t\in[0,\tau]}|x^\e_t(x)-x^0_t(x)|\geq \frac{\e^{2\g}}{2}\right) \,d\tau\\
&=\sum_{j=0}^{k_\e-1}\int_{j/\e^\theta}^{(j+1)/\e^\theta}\beta_\e e^{-\beta_\e \tau}
\P\left( \sup_{t\in[0,\tau]}|x^\e_t(x)-x^0_t(x)|\geq \frac{\e^{2\g}}{2}\right) \,d\tau\\
&\leq \sum_{j=0}^{k_\e-1}
\P\left( \sup_{t\in[0,(j+1)/\e^\theta]}|x^\e_t(x)-x^0_t(x)|
\geq \frac{\e^{2\g}}{2}\right) \int_{j/\e^\theta}^{(j+1)/\e^\theta}
\beta_\e e^{-\beta_\e \tau}\,d\tau\\
&\leq \sum_{j=0}^{k_\e-1}%\left( e^{-\beta_\e/\e^{\theta}}\right)^j\cdot
\P\left( \sup_{t\in[0,(j+1)/\e^\theta]}|x^\e_t(x)-x^0_t(x)|\geq \frac{\e^{2\g}}{2}\right)  .
\ea

For $j\geq 0$,
\ba
\label{eq:h}
&\P\left( \sup_{t\in[0,(j+1)/\e^\theta]}|x^\e_t(x)-x^0_t(x)|\geq \frac{\e^{2\g}}{2}\right) 
\leq
\P\left( \sup_{t\in[0,1/\e^\theta]}|x^\e_t(x)-x^0_t(x)|\geq \frac{\e^{2\g}}{5}\right) \\
&+\P\left( \sup_{t\in[0,1/\e^\theta]}|x^\e_t(x)-x^0_t(x)|< \frac{\e^{2\g}}{5}, 
\sup_{t\in[1/\e^\theta,2/\e^\theta ]}|x^\e_t(x)-x^0_{t-1/\e^\theta}(x^\e_{1/\e^\theta}(x))|\geq \frac{\e^{2\g}}{5}\right) \\
&+\cdots\\
&+ \P\Bigg( \sup_{t\in[(k-1)/\e^\theta,k/\e^\theta]}
|x^\e_t(x)-x^0_{t-(k-1)/\e^\theta}(x^\e_{(k-1)/\e^\g})|< \frac{\e^{2\g}}{5}, 
\text{ for } 0\leq k\leq j,\\
&\sup_{t\in[j/\e^\theta,(j+1)/\e^\theta ]}|x^\e_t(x)-x^0_{t-j/\e^\theta}(x^\e_{j/\e^\theta}(x))|
\geq \frac{\e^{2\g}}{5}\Bigg) \\
&\leq \P\left( \sup_{t\in[0,1/\e^\theta]}|x^\e_t(x)-x^0_t(x)|\geq \frac{\e^{2\g}}{5}\right) +
\sup_{y_1^-\leq y\leq y_1^+} \P\left( \sup_{t\in[0,1/\e^\theta ]}
|x^\e_t(y)-x^0_t(y))|\geq \frac{\e^{2\g}}{5}\right) +\cdots\\
&+\sup_{   y_j^-\leq y\leq y_j^+      } \P\left( \sup_{t\in[0,1/\e^\theta ]}
|x^\e_t(y)-x^0_t(y))|\geq \frac{\e^{2\g}}{5}\right) 
\ea
The sequence $y_k^\pm$ is determined by a recurrence formula. 
For small $\e$ and any $x\in[a+\e^\g,b-\e^\g]$, 
we know that $t_\e(x)+\hat{T}<1/\e^\theta$ and thus 
$|x^0_{1/\e^\theta}(x)|\leq \delta$. Moreover, $|x|e^{-m_2 t}\leq |x_t^0(x)|\leq |x|e^{-m_1 t}$, $t\geq 0$, $|x|\leq \delta$.
Define for $k\geq 2$
\ba
y_1^+&=x^0_{1/\e^\theta}(x)+\e^{2\g}/5, &&&   y_1^-&=x^0_{1/\e^\theta}(x)-\e^{2\g}/5,    \\
y_k^+&=
\begin{cases}
y_{k-1}^+ e^{-m_1/\e^\theta}+\e^{2\g}/5, \quad y_{k-1}^+\geq 0,\\
y_{k-1}^+ e^{-m_2/\e^\theta}+\e^{2\g}/5, \quad y_{k-1}^+\leq 0,
\end{cases} 
&&&
y_k^-&=
\begin{cases}
y_{k-1}^- e^{-m_2/\e^\theta}-\e^{2\g}/5, \quad y_{k-1}^-\geq 0,\\
y_{k-1}^- e^{-m_1/\e^\theta}-\e^{2\g}/5, \quad y_{k-1}^-\leq 0.
\end{cases} 
\ea
It is easy to see that for small $\e$ and $k\to\infty$,
$y^+_k\to \frac{\e^{2\g}/5}{1-e^{-m_1/\e^\theta}}\leq \frac{\e^{2\g}}{4}$,
$y^-_k\to -\frac{\e^{2\g}/5}{1-e^{-m_2/\e^\theta}}\geq -\frac{\e^{2\g}}{4}$.
%and the inequality \eqref{eq:h} holds.
Applying Proposition \ref{p:3} with $c=1/5$ and Lemma \ref{l:1} we get for some positive $p_1$ that
for $0<p\leq p_1$ and $\e\leq\e_2\leq\e_1$,
\be
\P\left( \sup_{t\in[0,(j+1)/\e^\theta]}|x^\e_t(x)-x^0_t(x)|\geq \frac{\e^{2\g}}{2}\right) 
\leq (j+1)\P\left( \sup_{t\in[0,1/\e^\theta]}|\e\xi^\e_t|\geq \e^{4\g}\right) 
\leq (j+1) e^{-1/\e^{p_1}}
\ee
and therefore
\ba
\int_0^{k_\e/\e^\theta}
&\beta_\e e^{-\beta_\e \tau}
\P\left( \sup_{t\in[0,\tau]}|x^\e_t(x)-x^0_t(x)|\geq \frac{\e^{2\g}}{2}\right) \,d\tau
\leq e^{-1/\e^{p_1}} \sum_{j=0}^{k_\e-1} (j+1) 
= \frac{1}{2}k_\e(k_\e+1)e^{-1/\e^{p_1}}.
\ea
Combining the latter formula with \eqref{eq:f} we obtain
the estimate needed for $0<p\leq p_0=\min\{\theta/2,p_1\}$, $\e\downarrow 0$.
\end{proof}

\subsection{Dynamics on unbounded interval, $a=-\infty$. Return from infinity\label{s:unb}}

In this section we show that with high probability the process $x^\e(x)$ reaches some 
fixed compact
neighbourhood of the origin in finite time. 

Recall that due to Assumption \textbf{U3} there is $N>0$ such that
$-U'(x)>c_1|x|^{1+c_2}$, for some $c_1$, $c_2>0$ and $x\leq -N$.  

Additionally, we assume that $N$ is sufficiently large, so that for any $x<-N$,
\be
 -|x|^{1+c_2}+|x+ \tfrac{1}{2}|^{1+c_2} +\tfrac{1}{4}(1+c_2)|x|^{c_2}<0.
\ee
Indeed this inequality holds, since for $x\to-\infty$,
\be
-|x|^{1+c_2}+|x+ \tfrac{1}{2}|^{1+c_2} +\tfrac{1}{4}(1+c_2)|x|^{c_2}=
-\tfrac{1}{4}(1+c_2)|x|^{c_2}+o(|x|^{c_2}).
\ee

We compare $x^\e(x)$ with the solution of the SDE
\be
v_t^\e(v)=v+c_1\int_0^t |v_{s-}^\e|^{1+c_2}\,ds+\e\xi_t^\e, \quad t\geq 0.
\ee

For some $M>N$ and $x\leq -M$ define stopping times
\ba
\tau_x=\inf\{t\geq 0: x^\e_t(x)\geq -M\},\\
\sigma_v=\inf\{t\geq 0: v^\e_t(v)\geq -M\}.
\ea

\begin{lemma}
For $v<x<-M$, $v_t^\e(v)<x_t^\e(x)$ a.s.\ on $t\in[0,\tau_x)$. 
\end{lemma}
\begin{proof}
Consider the difference
\be
\phi_t(x,v)=x_t^\e(x)-v_t^\e(v)=x-v+\int_0^t(-U'(x^\e_{s-})- c_1|v_{s-}^\e|^{1+c_2})\,ds.
\ee
The function $\phi_t(x,v)$
is absolutely continuous in $t$, $\phi_0(x,v)=x-y>0$. Let $t$ be the first time instant before
$\tau_x$ such that $\phi_0(x,v)=0$. This means that the left Dini derivative of $\phi$
is non-positive at $t$, $D^-\phi_t(x,v)=-U'(x^\e_{t-})- c_1|v_{t-}^\e|^{1+c_2}\leq 0$.
On the other hand, the processes $x^\e$ and $v^\e$ have the same jumps, so
$x^\e_{t}=v_{t}^\e$ if and only if $x^\e_{t-}=v_{t-}^\e$ which leads to a contradiction
with the assumptions.
\end{proof}

\begin{cor}
\label{c:comp}
For $v<x<-M$, $\tau_x\leq \sigma_v$ a.s.
\end{cor}

Fix some $M>N$ consider $T_M=\int_{-\infty}^{-M+1}\frac{dv}{c_1|v|^{1+c_2}}$.
Moreover, we can choose $M$ so that $v^0_{T_M}(-M)\leq -N$.

\begin{lemma}
\label{l:v-dev}
On the event 
$\mathcal{E}_{T_M}$ the following
holds a.s.\
\be
\sup_{t\in[0,T_M]}|v^\e_t(v)-v^0_t(v)|\leq 1
\ee
uniformly for $v\leq -M$.
\end{lemma}
\begin{proof}
As in Lemma \ref{l:z},
consider the representation $v^\e_t(v)=v^0_t(v)+\e w_t^\e(v)+r^\e_t(v)$ with 
\ba
v^0_t(v)&=v+c_1\int_0^t|v^0_s|^{1+c_2}\, ds,\\
w_t^\e(v)&=\xi^\e_t-c_1(1+c_2)|v^0_t(v)|^{1+c_2}\int_0^t\xi^\e_{s-}
\frac{ds}{|v^0_s(v)|}.
\ea
To estimate $w^\e$ we recall equations \eqref{eq:z1} and \eqref{eq:z2}
and immediately get
\be
\sup_{[0,T_M]}|w^\e_t(v)|\leq 2\sup_{[0,T_M]}|\xi^\e_t|, \quad v\leq -M.
\ee
The remainder term $r^\e$ satisfies the equation
\be
r_t^\e(v)=c_1\int_0^t\left(|v_s^0(v)+\e w_{s-}^\e(v)+ r_s^\e(v)|^{1+c_2}
-|v_s^0(v)|^{1+c_2}+(1+c_2)|v_s^0(v)|^{c_2} \e w_{s-}^\e(v)
 \right) \, ds
\ee
Assume, there exists a smallest $\tau\in[0,T_M]$ such that $r^\e_\tau(u)=3/4$.
Then the left Dini derivative of $r^\e$ at this point is non-negative, i.e.\
\be
|v_\tau^0(v)+\e w_{\tau-}^\e(v)+ \tfrac{3}{4}|^{1+c_2}
-|v_\tau^0(v)|^{1+c_2}+(1+c_2)|v_\tau^0(v)|^{c_2} \e w_{\tau-}^\e(v)\geq 0
\ee 
On the other hand on the event $\mathcal{E}_{T_M}$ we have  $|\e w_{s-}^\e(v)|< 1/4$ a.s.\ 
for $\e$ small enough, thus
\ba
&|v_\tau^0(v)+\e w_{\tau-}^\e(v)+ \tfrac{3}{4}|^{1+c_2}
-|v_\tau^0(v)|^{1+c_2}+(1+c_2)|v_\tau^0(v)|^{c_2} \e w_{\tau-}^\e(v)\\
&<|v_\tau^0(v)+ \tfrac{1}{2}|^{1+c_2}
-|v_\tau^0(v)|^{1+c_2}+\tfrac{1}{4}(1+c_2)|v_\tau^0(v)|^{c_2}<0,
\ea 
and a contradiction is reached. 

The estimate $r^\e_t\geq -3/4$ is obtained analogously, 
and the Lemma is proved.
\end{proof}

\begin{lemma}
\label{l:TK} 
For $x\leq-M$,
\be
\tau_x\leq T_M
\ee
a.s.\ on the event $\mathcal{E}_{T_M}$.
\end{lemma}
\begin{proof}
For any $x\leq -M$ compare $x^\e(x)$ with $v^\e(x-1)$. The the statement follows from
Corollary~\ref{c:comp},  Lemma~\ref{l:v-dev} and the definition of the time instant $T_M$.
\end{proof}

%\newpage

%%%%%%%%%%%%%%%%%%%%%%%%%%%%%%%%%%%%%%%%%%%%%%%%%%%%%%%%%%%%%%%%%%%%%%%%%%%%%%%%%%%%%%
%%%%%%%%%%%%%%%%%%%%%%%%%%%%%%%%%%%%%%%%%%%%%%%%%%%%%%%%%%%%%%%%%%%%%%%%%%%%%%%%%%%%%%
%%%%%%%%%%%%%%%%%%%%%%%%%%%%%%%%%%%%%%%%%%%%%%%%%%%%%%%%%%%%%%%%%%%%%%%%%%%%%%%%%%%%%%
%%%%%%%%%%%%%%%%%%%%%%%%%%%%%%%%%%%%%%%%%%%%%%%%%%%%%%%%%%%%%%%%%%%%%%%%%%%%%%%%%%%%%%
%%%%%%%%%%%%%%%%%%%%%%%%%%%%%%%%%%%%%%%%%%%%%%%%%%%%%%%%%%%%%%%%%%%%%%%%%%%%%%%%%%%%%%
%%%%%%%%%%%%%%%%%%%%%%%%%%%%%%%%%%%%%%%%%%%%%%%%%%%%%%%%%%%%%%%%%%%%%%%%%%%%%%%%%%%%%%
%%%%%%%%%%%%%%%%%%%%%%%%%%%%%%%%%%%%%%%%%%%%%%%%%%%%%%%%%%%%%%%%%%%%%%%%%%%%%%%%%%%%%%
%%%%%%%%%%%%%%%%%%%%%%%%%%%%%%%%%%%%%%%%%%%%%%%%%%%%%%%%%%%%%%%%%%%%%%%%%%%%%%%%%%%%%%
%%%%%%%%%%%%%%%%%%%%%%%%%%%%%%%%%%%%%%%%%%%%%%%%%%%%%%%%%%%%%%%%%%%%%%%%%%%%%%%%%%%%%%
%%%%%%%%%%%%%%%%%%%%%%%%%%%%%%%%%%%%%%%%%%%%%%%%%%%%%%%%%%%%%%%%%%%%%%%%%%%%%%%%%%%%%%

\subsection{Final estimate for $|x^\e_t-x^0_t|$, $a=-\infty$}

\begin{prop}
\label{p:12}
Let $a=-\infty$, $\rho\in(0,1)$.  Let $x^\e(x)$ and $x^0(x)$ satisfy \eqref{eq:eqs}. 
Let $T(\e)$ be an exponentially distributed random variable with mean $1/\beta_\e$
and let $\xi^\e$ and $T(\e)$ be independent. Then for any $\g\in (0,(1-\rho)/4)$.
there is $p_0>0$ and $\e_0>0$ such that the following estimate holds for all $0<\e\leq \e_0$
and $0<p\leq p_0$:
\ba
\sup_{x\leq b-\e^\g}
\Px\left( \sup_{t\in [0,\tau_x\wedge T(\e)]} x^\e_t(x)\geq -M+1\text{ or }
\sup_{t\in [\tau_x\wedge T(\e),T(\e)]}|x^\e_t(x)-x^0_{t-\tau_x\wedge T(\e)}(x^\e_{\tau_x}(x))|
\geq \frac{\e^{2\g}}{2}\right) 
\leq \exp{\left( -1/\e^p \right) }
\ea
where $\tau_x=\inf\{t\geq 0: x^\e_t(x)\geq -M\}$.
\end{prop}
\begin{proof}
For $x\in[-M,b-\e^\g]$ we have $\tau_x=0$ and the estimate
coincides with those of Proposition \ref{p:11} applied for a potential well 
$[-M,b-\e^\g]$, i.e.\ for the estimate holds for $0<\e\leq \e_1$ and $0< p\leq p_1$.

Consider the case $x\leq -M$. 
First due to Lemma \ref{l:1}, $\P(\mathcal{E}^c_{T_M})\leq e^{-1/\e^p}$, $0<\e\leq \e_2$, 
$0<p\leq p_2$.
Then, with the help of Markov property we obtain
\ba
\P&\left( \sup_{t\in [0,\tau_x\wedge T(\e)]} x^\e_t(x)\geq -M+1\text{ or }
\sup_{t\in [\tau_x\wedge T(\e),T(\e)]}|x^\e_t(x)-x^0_{t-\tau_x\wedge T(\e)}(x^\e_{\tau_x}(x))|
\geq \frac{\e^{2\g}}{2}\right)\\
&=\P(\mathcal{E}^c_{T_M})
+\P\left( \sup_{t\in [0,\tau_x\wedge T(\e)]} x^\e_t(x)\geq -M+1,\mathcal{E}_{T_M}\right)(=0) \\
&+\P\left(  \sup_{t\in [0,\tau_x\wedge T(\e)]} x^\e_t(x)<-M+1     ,   
\sup_{t\in [\tau_x\wedge T(\e),T(\e)]}|x^\e_t(x)-x^0_{t-\tau_x\wedge T(\e)}(x^\e_{\tau_x}(x))|
\geq \frac{\e^{2\g}}{2},\mathcal{E}_{T_M}\right) \\
&\leq \P(\mathcal{E}^c_{T_M})
+\sup_{y\in[-M,-M+1]}
\P\left( \sup_{t\in [0,T(\e)]}|x^\e_t(y)-x^0_t(y)|\geq \frac{\e^{2\g}}{2}\right) \leq e^{-1/\e^p}
\ea
for some positive $0<p\leq \min\{p_1,p_2\}$ and $0<\e\leq \min\{\e_1,\e_2\}$ small enough.
%%%
\end{proof}

\section{Exit from a single well\label{s:law:sigma}}

For $i=1,\dots,n$ consider the wells  of the potential 
$U$ with local minima at $m_i$. For $\e>0$ and $\g>0$ consider the following $\e$-dependent inner
neighbourhoods of the wells:
\ba 
\Omega^i&=(s_{i-1},s_i),\\
\Omega^i_\e&=[s_{i-1} +2\e^\g,s_i-2\e^\g],
\ea
where by convention $\Omega^1=(-\infty,s_1)$, $\Omega^1_\e=(-\infty,s_1-2\e^\g]$,
$\Omega^n=(s_{n-1},+\infty)$, and $\Omega^n_\e=[s_{n-1}+2\e^\g,+\infty)$.

Consider the following life times of the process $X^\e$ in the potential wells:
\ba
\sigma^i(\e)&=\inf\{t\geq 0: X^\e_t(\cdot)\notin [s_{i-1}+\e^\g, s_i-\e^\g]   \}, \quad i=1,\dots, n.
\ea

Let
\be
\lambda^i(\e)=H_-\left(\frac{s_{i-1}-m_i}{\e}\right) +H_+\left( \frac{s_i-m_i}{\e}\right), 
\quad i=1,\dots, n.
\ee

\begin{prop}
\label{p:q}
There exists $\g_0>0$ such that for any $0<\g\leq \g_0$, 
$x\in\Omega^i_\e$, $i=1,\dots,n$, any $C>0$ there exists $\e_0>0$ such that for $0<\e\leq \e_0$,
\ba
e^{-u(1+C)}(1-C) \leq \Px\left(   \lambda^i(\e)\sigma^i(\e)>u \right)\leq e^{-u(1-C)}(1+C).
\ea
Consequently,
\ba
\lim_{\e\downarrow 0}\lambda^i(\e) \E_x \sigma^i(\e)=1.
\ea
Moreover, for $j\neq i$,
\be
\label{eq:s3}
\lim_{\e\downarrow 0} \Px(X^\e_{\sigma^i(\e)}\in \Omega_\e^j)
%=\frac{\left||s_{j-1}-m_i|^{-r}-|s_j-m_i|^{-r}\right|}{|s_{i-1}-m_i|^{-r}+|s_i-m_i|^{-r}}
=\frac{q_{ij}}{q_i}.
\ee
% and
% \be
% \lim_{\e\downarrow 0}\lambda^i(\e) \E_x \left[ \sigma^i(\e)\I{X^\e_{\sigma^i(\e)}\in \Omega_\e^j} \right] 
% =\frac{\left||s_{j-1}-m_i|^{-r}-|s_j-m_i|^{-r}\right|}{|s_{i-1}-m_i|^{-r}+|s_i-m_i|^{-r}}.
% \ee
\end{prop}

Proposition \ref{p:q} will easily follow from
Lemmas~\ref{l:above} and  \ref{l:below} formulated below. The proof is rather 
technical and consists in applying the strong Markov property and accurate
estimations of certain probabilities.

\subsection{Useful technicalities \label{s:tech}}
\medskip

\subsubsection{Dynamics between big jumps}

Due to the strong Markov property, for any stopping time $\tau$
the process $\xi^\e_{t+\tau}-\xi^\e_\tau$, $t\geq 0,$ is
also a L\'evy process with the same law as $\xi^\e$.

For $k\geq 1$ consider processes
\begin{equation}
\begin{aligned}
\xi^k_t&=\xi^\e_{t+\tau_{k-1}}-\xi^\e_{\tau_{k-1}},\\
x^k_t(x)&=x-\int_0^t U'(x^k_{s-})\, ds +\e\xi^k_t, \quad t\in[0,T_k].
\end{aligned}
\end{equation}

In our notation, for $x\in\mathbb R$,
\begin{equation}
\begin{aligned}
X_t^{\e}&=x^1_t(x)+\e W_1^\e\I{t=T^\e_1}, \quad t\in[0,T^\e_1],\\
X_{t+\tau_1^\e}^{\e}&=x^2_t( x^1_{\tau_1^\e}+\e W_1)+\e W_2\I{t=T^\e_2}, \quad t\in[0,T^\e_2],\\
&\cdots\\
X_{t+\tau_{k-1}^\e}^{\e}&=x^k_t( x^{k-1}_{\tau_{k-1}^\e}+\e W^\e_{k-1})+\e W_k^\e\I{t=T_k^\e}, 
\quad t\in[0,T_k^\e].
\end{aligned}
\end{equation}

Denote $W_0^\e=T_0^\e=0$, $x^1(0)=x$, and write $\I{A}$ for the indicator function 
of a measurable set $A$.

\subsubsection{Constants $\rho$, $\g$ and $p_0$}

We assume that the threshold power $\rho$ and the constant $\g_0$ are fixed and satisfy 
\be
\label{eq:rhog}
\frac{1}{2}<\rho<1, \qquad 0<\g_0<\frac{1}{4}(1-\rho).
\ee
Then, for $0<\g\leq \g_0$ there is $p_0>0$ such that Propositions \ref{p:11} and \ref{p:12}
hold simultaneously for all wells $\Omega^i_\e$, $i=1,\dots, n$,  
for $0<p\leq p_0$ and $0<\e\leq \e_1$.
Further, we require that
\begin{itemize}
\item
$2\g<\rho<1-2\g$ \qquad (will be used in \textbf{Steps A1-2} and \textbf{A2-2} of Section \ref{S:Si:above} 
and \textbf{Steps B1-2} and \textbf{B2-2} of Section  \ref{S:Si:below}),\\
\item
$r(2\rho-1)+\g>0$ \qquad (will be used in \textbf{Step A2-2} of Section \ref{S:Si:above} 
and \textbf{Steps B1-2} and \textbf{B2-2} of Section \ref{S:Si:below}),
\end{itemize}
where $r>0$ is the index of regular variation of the tail of L\'evy measure (Assumption \textbf{L2}),
which obviously holds for $\rho$ and $\g$ satisfying \eqref{eq:rhog}

\subsubsection{Constant $c$}

Throughout this section we use a constant $c$ such that the following holds for
$\e\in (0,\e_0]$ for some $\e_0>0$:
\ba
\label{eq:c}
\sup_{y\in [s_{i-1}+\e^\g,s_i-\e^\g]}|X^0_t(y)-m_i|&\leq \tfrac{\e^{2\g}}{2}\text{ for } t\geq c|\!\ln\e|,
\quad i=1,\dots,n,\\
\sup_{|y-s_i|\geq \e^\g}|X^0_t(y)-s_i|&\geq \e^\g+2\e^{2\g}\text{ for }t\geq c\e^\g,
\quad i=1,\dots,n-1.
\ea
Let us show that these inequalities hold for some $c>0$. 
Let
$T(x,y)=\inf\{t\geq 0\,:\, X^0_t(x)=y\}$.
Then for any $i=1,\cdots,n$, and due to the properties of $U$ we need to show that
\ba
&T(s_{i-1}+\e^\g,m_i-\e^{2\g}/2),T(s_{i}-\e^\g, m_i+\e^{2\g}/2)\leq c|\!\ln\e|\quad i=1,\dots,n,\\
&T(s_i-\e^\g,s_i-\e^\g-2\e^{2\g}),T(s_i+\e^\g,s_i+\e^\g+2\e^{2\g})\leq c  \e^\g, \quad i=1,\dots,n-1.
\ea
what easily follows from nondegeneracy properties of potential's extremae (Assumption \textbf{U2}).

\subsubsection{Technical Lemma}

For definiteness, we assume as in Section \ref{s:onewell} that the well's minimum is located at
the origin, and denote well's boundaries as $-\infty\leq a<0<b<+\infty$. Denote
$\lambda(\e)=H_-\left(\frac{a}{\e}\right) +H_+\left( \frac{b}{\e}\right)$. Denote
\ba
I&=[a,b],\\
I_{\e,1}&=[a+\e^\g,b-\e^\g],\\
I_{\e,2}&=[a+\e^\g+\e^{2\g},b-\e^\g-\e^{2\g}].
\ea
if $a>-\infty$ and
\ba
I&=(-\infty,b],\\
I_{\e,1}&=(-\infty,b-\e^\g],\\
I_{\e,2}&=(-\infty,b-\e^\g-\e^{2\g}].
\ea
if $a=-\infty$.

For $y\in I_{\e,1}$, $j\geq 1$, we introduce the following events:
\ba
\label{eq:B}
A^j_y&=A^j(y)=\{x^j_s(y)\in I_{\e,1}, s\in [0,T_j), x^j_{T_j}(y)+\e W_j\in I_{\e,1} \},\\
A^{j,-}_y&=A^{j,-}(y)=
\{x^j_s(y)\in I_{\e,1}, s\in [0,T_j), x^j_{T_j}(y)+\e W_j\in I_{\e,2} \},\\
B^j_y&=B^j(y)=\{x^j_s(y)\in I_{\e,1}, s\in [0,T_j), x^j_{T_j}(y)+\e W_j\notin I_{\e,1} \},\\
E_y^j&=\{\omega:\sup_{t\in [0,T_j]}|x^j_t(y)-x^0_t(y)|\leq \frac{\e^{2\g}}{2}\}\quad \text{ if } a>-\infty,\\
E_y^j&=\{\omega:\sup_{t\in [0,\tau_y\wedge T_j]} x^j_t(y)\leq -M+1\text{ and }
\sup_{t\in [\tau_y\wedge T_j,T_j]}|x^j_t(y)-x^0_{t-\tau_y\wedge T_j}(x^j_{\tau_y}(y))|
\leq \frac{\e^{2\g}}{2}\}\quad \text{ if } a=-\infty.
\ea
with $M>0$ and $\tau_y$ defined in Section \ref{s:unb}.
Let also $A_y=A_y^1$, $A^-_y=A_y^{1,-}$, $B_y=B_y^1$, $E_y=E_y^1$.

Due to Propositions \ref{p:11} and \ref{p:12},  $\sup_{y\in I_{\e,1}} \P(E_y^c)\leq e^{-1/\e^p}$ for 
$0<p\leq p_0$ and $\e$ small enough. 
The following Lemma will also be used in the sequel.
\begin{lemma} 
\label{l:t1}
There exists a positive $\e_0$ such that the following holds true for all 
$0<\e\leq\e_0$ and $y\in I_{1,\e}$ 
\ba
1.&&&\I{A_y}\I{E_y}\I{T_1\geq c|\!\ln\e|}\leq \I{\e W_1\in I},\\
2.&&&\I{B_y}\I{E_y}\I{|\e W_1|>\tfrac{\e^{2\g}}{2}}\I{T_1\geq c|\!\ln\e|}
\leq \I{\e W_1\notin I_{\e,2}},\\
3.&&&\I{A^-_y}\I{E_y}\I{T_1\geq c|\!\ln\e|}%\\%\I{|\e W_1|\leq\tfrac{\e^{2\g}}{2}}\\
%&&&\qquad\qquad\qquad
\geq \I{E_y}\I{T_1\geq c|\!\ln\e|} %\I{|\e W_1|\leq\tfrac{\e^{2\g}}{2}}
\I{\e W_1\in [a+3\e^\g,b-3\e^\g]},\\
4.&&&\I{B_y}\I{E_y}\I{T_1\geq c|\!\ln\e|}
\geq \I{E_y}\I{T_1\geq c|\!\ln\e|}\I{\e W_1\notin [a-\e^\g-\e^{2\g},b+\e^\g+\e^{2\g}]}.
\ea
\end{lemma}
\begin{proof}
Essentially the statements follow from the fact that on 
$E_y\cap\{T_1\geq c|\!\ln\e|\}$, the inequality $|x^1_{T_1}(y)|\leq \e^{2\g}$ holds a.s.\ 
for all $y\in I_{\e,1}$. Indeed, if $a$ is finite this follows from
 Proposition~\ref{p:11}
and definition of the time $c|\!\ln\e|$. If $a=-\infty$, the statement follows from  
Proposition~\ref{p:12}.

\end{proof}

%%%
%%%
\subsection{Proof of Proposition \ref{p:q}. Upper estimate \label{S:Si:above}}
%%%
%%%

In this subsection we give an estimate of
$\Px(\lambda(\e)\sigma(\e)>u)$
%$\Px([\frac{\alpha\e^\alpha}{p^\alpha}]^{-1}\sigma(\e)>u)$
from above as $\e\to 0$, $u>0$. 

\begin{lemma}
\label{l:above}
For any $C>0$ there exist $\e_0>0$ such that for $0<\e\leq \e_0$ and $x\in I_{\e,2}$ 
\ba
\Px\left(\lambda(\e)\sigma(\e)>u\right)\leq e^{-u(1-C)}(1+C),
\ea
uniformly in $u\geq 0$.

\end{lemma}
\begin{proof}
For $x\in I_{\e,1}$, we use the following obvious inequality
\ba
\Px\left(\lambda(\e)\sigma(\e)>u\right)
&=\sum_{k=1}^\infty \P\left(\lambda(\e)\tau_k>u\right)
\Px(\sigma(\e)=\tau_k)\\
&+\Px\left(\lambda(\e)\sigma(\e)>u
\Big|\sigma(\e)\in (\tau_{k-1},\tau_k)\right)\Px(\sigma(\e)\in (\tau_{k-1},\tau_k))\\
&\leq  \sum_{k=1}^\infty
 \P\left(\lambda(\e)\tau_k>u\right)
\Big[ \Px(\sigma(\e)=\tau_k)+
\Px(\sigma(\e)\in (\tau_{k-1},\tau_k))\Big].
\ea

Then for any $x\in I_{\e,1}$ applying the independence 
and law properties of the
processes $x^j$, $j\in\mathbb N$, the following chain of inequalities is deduced which
results in a factorisation formula for the
probability under estimation: % (compare with \eqref{eq:fact}):
\begin{equation}
\label{eq:stauk}
\begin{aligned}
\Px(\sigma(\e)&=\tau_k)=\E_x\I{X_s^\e\in I_{\e,1}, s\in[0,\tau_k), X_{\tau_k}^\e\notin I_{\e,1}}\\
&= \E_x\prod_{j=1}^{k-1}
\I{A^j(X^\e_{\tau_{j-1}})} 
\cdot
\I{B^k(X^\e_{\tau_{k-1}})}
\leq \E\prod_{j=1}^{k-1}
\sup_{y\in I_{\e,1}}\I{A^j_y} 
\cdot
\sup_{y\in I_{\e,1}}\I{B^k_y}\\
&=\prod_{j=1}^{k-1}
\E\left[ \sup_{y\in I_{\e,1}}\I{A^j_y}\right]
\cdot
\E \left[ \sup_{y\in I_{\e,1}}\I{B^k_y}\right] 
=\left( \E\left[  \sup_{y\in I_{\e,1}}\I{A_y}\right]\right)^{k-1} 
\cdot
\E \left[ \sup_{y\in I_{\e,1}}\I{B_y}\right].
\end{aligned}
\end{equation}
Analogously we estimate the probability to exit between the $(k-1)$-th and the 
$k$-th arrival times of
the compound Poisson process $\eta^\e$, $k\in\mathbb N$.
Here we distinguish two cases.

In the first case, $k=1$, $x\in I_{\e,2}$. Then
\ba
&\Px(\sigma(\e)\in(\tau_0,\tau_1))=\Px(\sigma(\e)\in(0,T_1))=
\E_x\I{X_s^\e\notin I_{\e,1} \text{ for some } s\in (0,T_1)}\\
&\leq
\E\left[\sup_{y\in I_{\e,2} }\I{x^1_s(y)\notin I_{\e,2} \text{ for some } s\in[0,T_1]}\right].
\ea

In the second case, $k\geq 2$, $x\in I_{\e,1}$. Then
\ba
&\Px(\sigma(\e)\in(\tau_{k-1},\tau_k))=
\E_x\I{X_s^\e\in I_{\e,1}, s\in[0,\tau_{k-1}],
X_s^\e \notin I_{\e,1} \text{ for some } s\in (\tau_{k-1},\tau_k)}\\
&= \E_x \prod_{j=1}^{k-1}\I{A^j(X^\e_{\tau_{j-1}})} 
\cdot
\I{x^k_s(X^\e_{\tau_{k-1}})\notin I_{\e,1} \text{ for some } s\in [0,T_k)}\\
&\leq \E\prod_{j=1}^{k-2}
\sup_{y\in I_{\e,1}}\I{A^j_y}
\cdot
\sup_{y\in I_{\e,1}}\I{A^{k-1}_y}
\I{x^k_s(x^{k-1}_{T_{k-1}}(y)+\e W_{k-1})\notin I_{\e,1}\text{ for some } s\in[0,T_k]}\\
&=\left( \E\left[  \sup_{y\in I_{\e,1}}\I{A_y}\right]\right) ^{k-2}
\E \left[ \sup_{y\in I_{\e,1}}
\I{A_y}\I{x^2_s(x^1_{T_1}(y)+\e W_1)\notin I_{\e,1}\text{ for some } s\in[0,T_2]}\right].
\ea
Next we specify separately in four steps the further estimation for the four different
events appearing in the formulae for
$\Px(\sigma(\e)=\tau_k)$ and $\Px(\sigma(\e)\in(\tau_{k-1},\tau_k))$.

\noindent
\textbf{Step A1-1.} Consider  $\I{A_y}$. 
For $y\in I_{\e,1}$, we may estimate with help of Lemma~\ref{l:t1}
\ba
\label{eq:a1}
&\I{A_y}\leq \I{A_y}\I{E_y}+\I{E_y^c}
\leq \I{A_y}\I{E_y}\I{|\e W_1|>\tfrac{\e^{2\g}}{2}}
+\I{|\e W_1|\leq\tfrac{\e^{2\g}}{2}}+\I{E_y^c}\\
&\leq \I{A_y}\I{E_y}\I{|\e W_1|>\tfrac{\e^{2\g}}{2}}\I{T_1\geq c|\!\ln\e|}\\
&+\I{A_y}\I{E_y}\I{|\e W_1|>\tfrac{\e^{2\g}}{2}}\I{T_1< c|\!\ln\e|}
+\I{|\e W_1|\leq\tfrac{\e^{2\g}}{2}}+\I{E_y^c}\\
& \leq \I{|\e W_1|>\tfrac{\e^{2\g}}{2}}  \I{\e W_1 \in I}
+\I{|\e W_1|>\tfrac{\e^{2\g}}{2}} \I{T_1< c|\!\ln\e|}
+\I{|\e W_1|\leq\tfrac{\e^{2\g}}{2}}+ \I{E_y^c}\\
&= \I{\e W_1 \in I}
+\I{|\e W_1|>\tfrac{\e^{2\g}}{2}}\I{T_1< c|\!\ln\e|}
+ \I{E_y^c}.
\ea

\noindent
\textbf{Step A2-1.} Consider  $\I{B_y}$.
For $y\in I_\e^1$, we may estimate with help of Lemma~\ref{l:t1} 
\ba\label{eq:a2}
&\I{B_y}\leq \I{B_y}\I{E_y}+\I{E_y^c} \\
&=\I{B_y}\I{E_y}\I{|\e W_1|>\tfrac{\e^{2\g}}{2}}\I{T_1\geq c|\!\ln\e|}
+\I{B_y}\I{E_y}\I{|\e W_1|>\tfrac{\e^{2\g}}{2}}\I{T_1< c|\!\ln\e|}\\
&+\I{B_y}\I{E_y}\I{|\e W_1|\leq\tfrac{\e^{2\g}}{2}}\I{T_1\geq c\e^\gamma}
+\I{B_y}\I{E_y}\I{|\e W_1|\leq\tfrac{\e^{2\g}}{2}}\I{T_1 < c\e^\gamma}
+\I{E_y^c} \\
&\leq\I{\e W_1\notin I_{\e,2}}+\I{|\e W_1|>\tfrac{\e^{2\g}}{2}}\I{T_1< c|\!\ln\e|} 
+0+ \I{T_1< c\e^\gamma} +\I{E_y^c}.
\ea

\noindent
\textbf{Step A3-1.} Consider  $\I{x^1_s(y)\notin I_{\e,1} \text{ for some } s\in[0,T_1]}$.
For $y\in I_{\e,2}$, we may estimate
\ba
\label{eq:a3}
&\I{x^1_s(y)\notin I_{\e,1}\text{ for some } s\in[0,T_1]}\leq
\I{E_y^c}+
\I{x^1_s(y)\notin I_{\e,1} \text{ for some } s\in[0,T_1]} \I{E_y}
=\I{E_y^c}
\ea

\noindent
\textbf{Step A4-1.} Consider 
$\I{A_y}
\I{x^2_s(x^1_{T_1}(y)+\e W_1)\notin I_{\e,1}\text{ for some } s\in[0,T_2]}$
for $y\in I_{\e,1}$, we may estimate
\ba
\label{eq:a4}
&\I{A_y}
\I{x^2_s(x^1_{T_1}(y)+\e W_1)\notin I_{\e,1}\text{ for some } s\in[0,T_2]}\\
&=\I{x^1_s(y)\in I, s\in(0,T_1], x^1_{T_1}(y)+\e W_1\leq \in I_{\e,2},
x^2_s(x^1_{T_1}(y)+\e W_1)\in I_{\e,1}\text{ for some } s\in[0,T_2]}\\
&+\I{x^1_s(y)\in I, s\in(0,T_1], x^1_{T_1}(y)+\e W_1\in I_{\e,1}\backslash I_{\e,2},
x^2_s(x^1_{T_1}(y)+\e W_1)\notin I_{\e,1}\text{ for some } s\in[0,T_2]}\\
&\leq \I{x^1_s(y)\in I, s\in(0,T_1], x^1_{T_1}(y)+\e W_1\in I_{\e,2} }\cdot
\sup_{y\in I_{\e,2}}\I{x^2_s(y)\notin I_{\e,1}\text{ for some } s\in[0,T_2]}\\
&+\I{x^1_s(y)\in I, s\in(0,T_1], x^1_{T_1}(y)+\e W_1\in I_{\e,1}\backslash I_{\e,2}  }\\
&\leq \sup_{y\in I_{\e,2}}\I{x^2_s(y)\notin  I_{\e,1} \text{ for some } s\in[0,T_2]}\\
&+\I{x^1_s(y)\in I_{\e,1}, s\in(0,T_1], x^1_{T_1}(y)+\e W_1\in I_{\e,1}\backslash I_{\e,2}}.
\ea
The first term in the resulting expression in the \textbf{Step A4-1} is 
identical to the expression handled in \textbf{Step A3-1},
while the second term requires an inessential modification of the estimation in \textbf{Step A2-1}, namely
we consider an event $\{x^1_{T_1}(y)+\e W_1\in I_{\e,1}\backslash I_{\e,2}\}$ instead of
$\{x^1_{T_1}(y)+\e W_1\notin I_{\e,1}\}$.

Now we apply \eqref{eq:a1}, \eqref{eq:a2}, \eqref{eq:a3} and \eqref{eq:a4} to
estimate the expectations
treated in \textbf{Steps A1-1} --- \textbf{A1-4} above.
Let $C$ be a positive constant.
%Fix also $0<\delta<\min\{\gamma/2,\alpha(1-\gamma),\alpha/2\}=\min\{\gamma/2,\alpha/2\}$.

\noindent\textbf{Step A1-2.} Estimate
$\E\left[  \sup_{y\in I_{\e,1}}\I{A_y}\right]$. We get for $2\g<\rho<1-2\g$, some $\e_1>0$ and all $\e\leq \e_1$ that
\ba
\label{eq:esta1}
\E\left[  \sup_{y\in I_{\e,1}}\I{A_y}\right]&\leq 
\P(\e W_1 \in I)+\P(|\e W_1|>\tfrac{\e^{2\g}}{2})\P(T_1< c|\!\ln\e|)
+\sup_{y\in I_{\e,1}}\P({E_y^c})\\
&1-\frac{H_-(a/\e)+H_+(b/\e)}{\beta_\e}
+\frac{H(1/(2\e^{1-2\g}))}{\beta_\e}\cdot c \beta_\e|\!\ln\e|+e^{-1/\e^p}\\
&\leq 1- \frac{H_-(a/\e)+H_+(b/\e)}{\beta_\e} 
\left(1-  \frac{c \beta_\e H(1/(2\e^{1-2\g}))|\!\ln\e|  +\beta_\e e^{-1/\e^p}}{H_-(a/\e)+H_+(b/\e)} \right)\\
&\leq 1- \frac{\lambda(\e)}{\beta_\e} \left(1-  \frac{C}{5}\right).
\ea

\noindent\textbf{Step A2-2.} Estimate
$\E \left[ \sup_{y\leq I_{\e,1}}\I{B_y}\right]$. In fact, for $r(2\rho-1)+\g>0$ and $2\g<\rho<1-2\g$ and $\e\leq \e_2$
\ba
\label{eq:esta2}
\E &\left[ \sup_{y\leq I_{\e,1}}\I{B_y}\right]
\leq \P(\e W_1 \notin I_{\e,2})
+\P(|\e W_1|>\tfrac{\e^{2\g}}{2})\P(T_1< c|\!\ln\e|)+\P(T_1<c\e^\g)
+ \sup_{y\leq I_\e^1}\P({E_y^c})\\
&\leq
\frac{H_-((a+\e^\g+\e^{2\g})/\e)+H_+((b-\e^\g-\e^{2\g})/\e)}{\beta_\e}
+ c H(1/(2\e^{1-2\g})) |\!\ln\e|
+c \beta_\e \e^\g
+e^{-1/\e^p}\\
&= \frac{H_-(a/\e)+H_+(b/\e)}{\beta_\e }  \\
&\times\left( \frac{H_-((a+\e^\g+\e^{2\g})/\e)+H_+((b-\e^\g-\e^{2\g})/\e)}{H_-(a/\e)+H_+(b/\e) }     
+  \frac{c \beta_\e H(1/(2\e^{1-2\g}))|\!\ln\e| 
+c\beta_\e^2\e^\g +\beta_\e e^{-1/\e^p}}{H_-(a/\e)+H_+(b/\e)}    \right)\\
&\leq 
\frac{\lambda(\e)}{\beta_\e }\left(1+\frac{C}{5}\right).
\ea
On this step to estimate the ratio $H_+((b-\e^\g-\e^{2\g})/\e)/H_+(b/\e)$ 
we used the uniform convergence of slowly varying functions, see Proposition \ref{p:usv}.

\noindent\textbf{Step A3-3.} Estimate
$\E
\left[ \sup_{y\in I_{\e,2}}\I{x^1_s(y)\notin I_{\e,1}\text{ for some } s\in[0,T_1]}\right]$. 
We have for $\e\leq \e_3$
\ba
\label{eq:esta3}
\E \left[ \sup_{y\in I_{\e,2}}\I{x^1_s(y)\notin I_{\e,1}\text{ for some } s\in[0,T_1]}\right]&
\leq \sup_{y\in I_{\e,2}}\P({E_y^c})\leq e^{-1/\e^p}
\leq \frac{C}{5}\cdot \frac{\lambda(\e)}{\beta_\e }.
\ea

\noindent\textbf{Step A4-2.} Estimate 
 $\E\left[ \sup_{y\in I_{\e,1}}\I{A_y}
\I{x^2_s(x^1_{T_1}(y)+\e W_1)\notin I_{\e,1}\text{ for some } s\in[0,T_2]}\right] $.
We finally obtain for $\e\leq \e_4$
\ba
\label{eq:esta4}
&\E\left[\sup_{y\in I_{\e,1}}\I{A_y}
\I{x^2_s(x^1_{T_1}(y)+\e W_1)\notin I_{\e,1}\text{ for some } s\in[0,T_2]}\right]\\
&\leq
\P({E_y^c})
+\P(\e W_1 \in [a+\e^\g-\e^{2\g},a+\e^\g+2\e^{2\g}])+\P(\e W_1 \in [b-\e^\g-2\e^{2\g},b-\e^\g+\e^{2\g}])\\
&+ \P(|\e W_1|>\tfrac{\e^{2\g}}{2})\P(T_1< c|\!\ln\e|)+ \sup_{y\in I_{\e,1}}\P({E_y^c})
\leq \frac{C}{5}\cdot \frac{\lambda(\e)}{\beta_\e }.
\ea
%%%
%%%
%%%%%%
%%%
Then for $x\in I_{\e,2}$, $0<\e\leq \e_0$ and $\e\leq \e_5<\min\{\e_1,\e_2,\e_3,\e_4\}$,
\ba
\Px&\left(\lambda(\e) \sigma(\e)>u\right) \leq
\P\left(\lambda(\e) \tau_1>u\right)  \frac{\lambda(\e)}{\beta_\e }
\left(1+\frac{2C}{5}\right) \\
&+\sum_{k=2}^\infty \P\left( \lambda(\e)\tau_k>u\right) 
\left[  1- \frac{\lambda(\e)}{\beta_\e} \left(1-  \frac{C}{5}\right)             \right]^{k-1}
\frac{\lambda(\e)}{\beta_\e}
\left[1+\frac{C}{5}
+\frac{C/5}{ 1- \frac{\lambda(\e)}{\beta_\e} \left( 1-  \frac{C}{5}\right)       }
 \right] \\
&\leq \sum_{k=1}^\infty
\int_u^\infty \frac{\beta_\e}{\lambda(\e)} 
e^{-\tfrac{\beta_\e}{\lambda(\e)}  t}
\frac{\left(  \tfrac{\beta_\e}{\lambda(\e)} t\right)^{k-1}}{(k-1)!}\, dt
\left[ 1- \frac{\lambda(\e)}{\beta_\e} \left(1-  \frac{C}{5}\right)            \right]^{k-1}
\frac{\lambda(\e)}{\beta_\e} \left(1+  \frac{3C}{5}\right)      \\
&\leq \left(1+\frac{3C}{5}\right)
\int_u^\infty  e^{- t (1-C/5)}\,dt
\leq \frac{1+3C/5}{1-C/5}e^{- u (1-C/5)}
\leq e^{- u (1-C)}(1+C).
\ea

In the previous formula we have changed summation and 
integration. This can be done
due to uniform convergence of the series which follows from dominated convergence.

% In the previous formula we have changed summation and integration. This can be done
% due to the uniform convergence of the series $\sum_{k=1}^\infty
% \beta_\e e^{-\beta_\e t}\frac{(\beta_\e t)^{k-1}}{(k-1)!}
% \left[ 1-\e^{\alpha/2}\frac{\theta}{2}(1-C_5\e^\delta)\right]^{k-1}$ for all $t\geq 0$ and
% $\e\leq \e_0$. Indeed, let $t^*_k$ be the coordinate of the maximum of the density
% of  the $(\beta_\e, k)$-Gamma distribution. For $k\geq 2$, it is easy to see that
% $t^*_k=\frac{k-1}{\beta_\e}$
% Then, with help of Stirling's formula we obtain, that
% \ba
% 0\leq \beta_\e e^{-\beta_\e t}\frac{(\beta_\e t)^{k-1}}{(k-1)!}\leq
% \beta_\e e^{-(k-1)}\frac{(k-1)^{k-1}}{(k-1)!}
% \sim \frac{1}{\sqrt{2\pi}}\frac{\beta_\e}{\sqrt{k-1}},\quad k\to \infty.
% \ea
% Then, for all $\e\leq \e_0$,
% \ba
% \sum_{k=1}^\infty
% \beta_\e e^{-\beta_\e t}\frac{(\beta_\e t)^{k-1}}{(k-1)!}
% \left[ 1-\e^{\alpha/2}\frac{\theta}{2}(1-C_5\e^\delta)\right]^{k-1}\leq
% c_1\frac{\beta_\e}{\e^{\alpha/2}\frac{\theta}{2}(1-C_5\e^\delta)}\leq c,
% \ea
% where the constant $c$  does not depend on $t$ and $\e$. Thus uniform convergence
% follows from dominated convergence.

\end{proof}

\subsection{Proof of Proposition \ref{p:q}. Lower estimate \label{S:Si:below}}

In this subsection we estimate $\Px(\lambda(\e)\sigma(\e)>u)$  
from below as $\e\to 0$, $u>0$. This leads to
the following Lemma with a rather technical proof again.

\begin{lemma}
\label{l:below}
For any $C>0$ there exist $\e_0>0$ such that
for all $0<\e\leq \e_0$ and  
$x\in I_{\e,2}$ 
\ba
\Px&\left(\lambda(\e)\sigma(\e)>u\right) \geq
e^{-u(1+C)}(1-C)
\ea
uniformly in $u\geq 0$.
\end{lemma}

\begin{proof}
We use the following inequality:
\ba
\Px\left(\lambda(\e)\sigma(\e)>u\right) 
% &=\sum_{k=1}^\infty \P(\lambda(\e)\tau_k>u)\Px(\sigma(\e)=\tau_k)+
% \Px\left(\lambda(\e)\sigma(\e)>u\Big|\sigma(\e)\in (\tau_{k-1},\tau_k)\right) 
% \Px(\sigma(\e)\in (\tau_{k-1},\tau_k))\\
\geq  \sum_{k=1}^\infty \P\left( \lambda(\e)\tau_k>u\right) 
 \Px(\sigma(\e)=\tau_k).
\ea

With arguments analogous to \eqref{eq:stauk} we obtain the factorization
\ba
&\Px(\sigma(\e)=\tau_k)=\E_x\I{X_s^\e\in I_{\e,1}, s\in[0,\tau_k), X_{\tau_k}^\e\notin I_{\e,1}}\\
&= \E_x \prod_{j=1}^{k-1}
\I{A^j(X^\e_{\tau_{j-1}})}\cdot\I{B^k(X^\e_{\tau_{k-1}})}\\
&\geq \E_x \prod_{j=1}^{k-1}
\I{A^{j,-}(X^\e_{\tau_{j-1}})}\cdot\I{B^k(X^\e_{\tau_{k-1}})}\\
&=\E\left[\prod_{j=1}^{k-1}
 \inf_{y \in I_{\e,2}  }\I{A^{j,-}_y} \cdot
\inf_{y \in I_{\e,2} }\I{B^k_y}\right] \\
&=\left( \E\left[  \inf_{y\in I_{\e,2}}\I{A^-_y} \right]\right)^{k-1} \cdot
\E \left[  \inf_{y\in I_{\e,2}}\I{B_y}   \right].
\ea
For $y\in I_{\e,2}$, 
we next specify separately in two steps the further estimation for the two different events 
appearing in the formulae for $\Px(\sigma(\e)=\tau_k)$.

\noindent
\textbf{Step B1-1.}
Consider the event $\I{A_y^-}$. We may estimate with help of Lemma~\ref{l:t1} 
\ba
\label{eq:b1}
&\I{A_y^-}\geq \I{A_y^-}\I{E_y}\\
&\geq
\I{A_y^-}\I{E_y}\I{|\e W_1|\leq\tfrac{\e^{2\g}}{2}}\I{T_1\geq c\e^\g}
+\I{A_y^-}\I{E_y}\I{|\e W_1|>\tfrac{\e^{2\g}}{2}}\I{T_1\geq c|\!\ln\e|}\\
&\geq
\I{E_y}\I{|\e W_1|\leq\tfrac{\e^{2\g}}{2}}\I{T_1\geq c\e^\g}
+\I{E_y}\I{|\e W_1|>\tfrac{\e^{2\g}}{2}}\I{T_1\geq c|\!\ln\e|}
\I{\e W_1\in [a+3\e^\g,b-3\e^\g]}\\
&\geq
\I{|\e W_1|\leq\tfrac{\e^{2\g}}{2}}\I{T_1\geq c\e^\gamma}
+\I{|\e W_1|>\tfrac{\e^{2\g}}{2}}\I{T_1\geq c|\!\ln\e|}
\I{\e W_1\in [a+3\e^\g,b-3\e^\g]}
-2\I{E_y^c}\\
&\geq\I{|\e W_1|\leq\tfrac{\e^{2\g}}{2}}
+\I{|\e W_1|>\tfrac{\e^{2\g}}{2}} \I{\e W_1 \in [a+3\e^\g,b-3\e^\g]}\\
&-\I{ T_1< c\e^\gamma}
-\I{|\e W_1|>\tfrac{\e^{2\g}}{2}}\I{ T_1< c|\!\ln\e|}-2\I{E_y^c}\\
&=\I{\e W_1 \in[a+3\e^\g,b-3\e^\g]}
-\I{ T_1< c\e^\gamma}
-\I{|\e W_1|>\tfrac{\e^{2\g}}{2}}\I{ T_1< c|\!\ln\e|}-2\I{E_y^c}
\ea

\noindent
\textbf{Step B2-1.}
With help of Lemma~\ref{l:t1} the event $\I{B_y}$ may be estimated as follows 
\ba
\label{eq:b2}
\I{B_y} &\geq \I{B_y} \I{E_y}\I{T_1\geq c|\!\ln\e|}\\
&\geq \I{E_y}\I{T_1\geq c|\!\ln\e|}\I{\e W_1\notin [a-\e^\g-\e^{2\g},b+\e^\g+\e^{2\g}]}\\
&\geq \I{\e W_1 \notin [a-\e^\g-\e^{2\g},b+\e^\g+\e^{2\g}] }\left(1 -\I{T_1< c|\!\ln\e|}
-\I{E_y^c  }\right).
\ea
Now we apply \eqref{eq:b1} and \eqref{eq:b2}  to estimate the expectations
appearing in the formula for $\Px(\sigma(\e)=\tau_k)$.
Let $C>0$.

\noindent\textbf{Step B1-2.} Here we estimate
$\E\left[ \inf_{y\in I_{\e,2} }\I{A_y^-}\right]$, $2\g<\rho<1-2\g$, $r(2\rho-1)+\g>0$. There exists
$\e_1>0$ such that for $0<\e\leq \e_1$ the following holds similarly to \eqref{eq:esta1} and  \eqref{eq:esta2}
\ba
&\E\left[ \inf_{y \in I_{\e,2}  }\I{A_y^-}\right]\\
&\geq  \P(\e W_1\in[a+3\e^\g,b-3\e^\g])
-\P(T_1< c\e^\gamma)-
\P(|\e W_1|>\tfrac{\e^{2\g}}{2})\P(T_1< c|\!\ln\e|)-2\sup_{y\in I_{\e,2} }\P(E_y^c)\\
&\geq 1-\frac{H_-((a+3\e^\g)/\e)+H_+((b-3\e^\g)/\e) }{\beta_\e}
-c \beta_\e \e^\g
-c H(1/(2\e^{1-2\g}))|\!\ln\e|- 2e^{-1/\e^p}\\
% &= 1- \frac{H_-(a/\e)+H_+(b/\e)}{\beta_\e}\\
% &\times
% \left( 
% \frac{H_-((a+\e^\g+\e^{2\g})/\e)+H_+((b-\e^\g-\e^{2\g})/\e)}{H_-(a/\e)+H_+(b/\e)}
% +\frac{  c \beta_\e^2 \e^\g
% +c \beta_\e H(1/(2\e^{1-2\g}))|\!\ln\e|+2\beta_\e e^{-1/\e^p}  }{H_-(a/\e)+H_+(b/\e)}
% \right) \\
&\geq 1- \frac{\lambda(\e)}{\beta_\e}\left(1+\frac{C}{2} \right).
\ea
Here we again used the uniform convergence from Proposition \ref{p:usv}. 

\noindent\textbf{Step B2-2.} We next estimate
$\E\left[\inf_{y\in I_{\e,2}}\I{B_y}\right]$, for which we obtain similarly for $0<\e\leq \e_2$ with 
some $\e_2>0$. 
\ba
\E\left[\inf_{y\in I_{\e,2}}\I{B_y}\right] &\geq \P(\e W_1\notin [a-\e^\g-\e^{2\g},b+\e^\g+\e^{2\g} ])
\left(1 -\P(T_1< c|\!\ln\e|)-\sup_{y\in I_{\e,2} }\P(E_y^c)\right)\\
% &\geq \frac{H_-(a/\e)+H_+(b/\e)}{\beta_\e}\\
% &\times\left(1-\frac{H_-(a/\e)-H_-((a-\e^\g-\e^{2\g})/\e) +H_+(b/\e)-H_+((b+\e^\g+\e^{2\g})/\e) }
% {H_+(a/\e)+H_+(b/\e)}\right)
% \left( 1-\frac{C}{4}\right)  \\
&\geq \frac{\lambda(\e)}{\beta_\e}\left( 1-\frac{C}{4}\right)^{\!2}\geq 
\frac{\lambda(\e)}{\beta_\e}\left( 1-\frac{C}{2}\right) .
\ea
Consequently for  $0<\e\leq \min\{\e_1,\e_2\}$ and   $x\in I_{\e,2}$,  
\ba
\Px\left( \lambda(\e)\sigma(\e)>u\right) 
&\geq \sum_{k=1}^\infty
\int_u^\infty \frac{\beta_\e}{\lambda(\e)} 
e^{- \tfrac{\beta_\e}{\lambda(\e)} t}
\frac{\left(\tfrac{\beta_\e}{\lambda(\e)}  t\right)^{k-1}}{(k-1)!}\, dt
\left[ 1-\frac{\lambda(\e)}{\beta_\e}\left(1+\frac{C}{2}\right)\right]^{k-1}
\frac{\lambda(\e)}{\beta_\e}\left(1-\frac{C}{2}\right)\\
&\geq \left( 1-\frac{C}{2}\right) 
\int_u^\infty  e^{- t (1+C/2)}\,dt
\geq \frac{1-C/2}{1+C/2}e^{- u (1+C/2)}
\geq e^{- u (1+C)}(1-C).
\ea
See the end of the proof of Lemma \ref{l:above} 
for the justification of switching the order
of summation and integration in the above argument. 
\end{proof}

\bigskip

\begin{Aproof}{Proposition}{\ref{p:q}}
The first statement of Proposition \ref{p:q} follows directly from Lemmas \ref{l:above}
and \ref{l:below}.

The estimate for the expected value of $\sigma_x^i(\e)$ follows easily from the equality
\be
\lambda^i(\e)\E_x \sigma^i(\e)=
\int_0^\infty \Px\left( \lambda^i(\e)\sigma^i(\e)>u\right)\,du.
\ee
To obtain the 
third statement we repeat the steps of the argument of Lemmas \ref{l:above}
and \ref{l:below} taking $u=0$ and redefining the event $B^j_y$ in \eqref{eq:B} and thereafter as
\be
\{x^j_s(y)\in I_{\e,1}, s\in [0,T_j], x^j_{T_j}(y)+\e W_j\in \Omega_\e^j\}.
\ee
Then, it is easy to see that for $x\in \Omega^i_\e$
\ba
&\left[ \frac{H_-(\frac{s_j-m_i}{\e})-H_-(\frac{s_{j-1}-m_i}{\e})}
{H_-(\frac{s_{i-1}-m_i}{\e})+H_+(\frac{s_j-m_i}{\e})}\right]^{-1} 
\Px(X^\e_{\sigma^i(\e)}\in \Omega^j_\e)\to 1, \quad \text{if } j<i,\\
&\left[\frac{H_+(\frac{s_{j-1}-m_i}{\e})- H_+(\frac{s_j-m_i}{\e})}
{H_+(\frac{s_{i-1}-m_i}{\e})+H_+(\frac{s_j-m_i}{\e})}\right]^{-1} 
\Px(X^\e_{\sigma^i(\e)}\in \Omega^j_\e)\to 1, \quad \text{if } i<j.
\ea
and the ratios in brackets converge to $q_{ij}/q_i$ as defined in \eqref{eq:q}.

\end{Aproof}

\section{Transitions between the wells\label{s:transitions}}

For $0<\Delta<\Delta_0=\min_{1\leq i\leq n}\{|m_i-s_{i-1}|,|m_i-s_i|\}$ and $x\in\mathbb{R}$ denote $B_\Delta(x)=\{y:|x-y|\leq \Delta\}$.
Consider the following stopping times:
\begin{align}
T^i(\e)&=
\inf\{t\geq 0: X^\e_t(\cdot)\in \bigcup_{k\neq i}\Omega^k_\e\},\\
\tau^i(\e)&=
\inf\{t\geq 0: X^\e_t(\cdot)\in \bigcup_{k\neq i} B_\Delta(m_k)  \}, \\
S^i(\e)&= \inf\{t\geq 0: X^\e_t(\cdot)\notin B_{2\e^\g}(s_i) \}, \quad i=1,\dots, n-1.
\end{align}
For $x\in \Omega^i_\e$, $T^i$ is the transition time between the wells. For $x\in B_\Delta(m_i)$, $\tau^i$ is the transition time between $\Delta$-neighbourhoods of
wells' minima, and for $x\in B_{2\e^\g}(s_i)$, $S_x^i$ is the exit time from a neighbourhood of the saddle point.

\begin{lemma}
\label{l:s}
Let $i=1,\dots,n-1$ and $x\in B_{2\e^\g}(s_i)$. Then
\ba
\lim_{\e\downarrow 0}  H(1/\e) \E_x S^i(\e) =0,\\
\ea
\end{lemma}
\begin{proof}
To estimate $\E_x S^i(\e)$ we notice that for $x\in B_{2\e^\g}(s_i)$, 
\be
S_x^i(\e)\leq \inf\{t>0: |\e L_t-\e L_{t-}|>4\e^\g\}=J(\e)
\text{ a.s.},
\ee
i.e.\ the first exit time of $X^\e$ from the $2\e^\g$-neighbourhood of the saddle point $s_i$
is a.s.\ bounded from above by the time of the first jump of $\e L$ exceeding $4\e^\g$.
Note that $J(\e)$ is exponentially distributed with mean
\be
\label{eq:S1}
\E J(\e)=\left( \int_{|y|>4/\e^{1-\g}}\nu(dy)\right)^{\!\!-1}=\frac{1}{H(4/\e^{1-\g})}.
\ee
The statement of the Lemma follows from the fact that $H(1/\e)/H(4/\e^{1-\g})\to 0$ as $\e\downarrow 0$. 
\end{proof}

\bigskip

\begin{prop}
\label{p:gT}
For $x\in\Omega^i_\e$ and $j\neq i$
\begin{align}
\label{eq:T1}
&\lim_{\e\downarrow 0}\Px(X_{T^i(\e)}^\e\in \Omega_\e^j) =\frac{q_{ij}}{q_i}\\
\label{eq:T2}
&\lim_{\e\downarrow 0}\Px(T^i(\e)>\sigma^i(\e)) =0,\\
\label{eq:T3}
& \lim_{\e\downarrow 0}\lambda^i(\e) \E_x T^i(\e)=1.
\end{align}
\end{prop}
\begin{proof}
It is obvious that for all $x\in\Omega_\e^i$
\be
\sigma^i(\e)\leq T^i(\e)\quad \Px\text{-a.s.}
\ee
We have the inequality 
\ba 
\Px(X_{T^i(\e)}^\e\in \Omega_\e^j)= \Px(X_{\sigma^i(\e)}^\e\in \Omega_\e^j)+
\Px(X_{T^i(\e)}^\e\in \Omega_\e^j, T^i(\e)>\sigma^i(\e))  \geq \Px(X_{\sigma^i(\e)}^\e\in \Omega_\e^j) .
\ea
Recall \eqref{eq:s3} in Proposition \ref{p:q} and
note that $\sum_{j\neq i}\frac{q_{ij}}{q_i}=1$. Then the limits \eqref{eq:T1} and \eqref{eq:T2} follow. 

For any $\delta>0$ there exists $\e_0>0$ such that for $0<\e\leq \e_0$ 
the following estimates hold
\ba
&\sup_{x\in \Omega^i_\e}\Px\left(X^\e_{\sigma^i(\e)}\in \bigcup_{j=1}^{n-1} B_{2\e^\g}(s_j)\right)
\leq \delta, \\
&\sup_{x\in \Omega^i_\e}\lambda^i(\e)\E_x\sigma^i(\e)\leq 1+\delta,\\
&\max_{1\leq j\leq n-1}\sup_{x\in B_{2\e^\g}(s_j)}\lambda^i(\e)\E_x S^j(\e)\leq \delta.
\ea
Then is easy to see  that
\ba
\lambda^i(\e)\E_x T^i(\e)
&\leq  \lambda^i(\e)\E_x \sigma^i(\e)+ \sum_{k=1}^\infty (k+1)(1+\delta)+k \delta)\delta^k
\leq 1+ \delta\cdot\Const
\ea
which proves \eqref{eq:T3}.
\end{proof}

\bigskip

\begin{prop}
\label{p:1f}
For any  $0<\Delta<\Delta_0$ the following limits hold
\begin{align}
\label{eq:t1}
&\lim_{\e\downarrow 0}\Px\left( X_{\tau^i(\e)}^\e\in B_\Delta(m_j)\right)  =\frac{q_{ij}}{q_i}\\
\label{eq:t2}
&\lambda^i(\e)\tau^i(\e)\stackrel{\mathcal{D}}{\to} \exp(1),\\
\label{eq:t3}
& \lim_{\e\downarrow 0}\lambda^i(\e) \E_x \tau^i(\e)=1.
\end{align}
uniformly for $x\in B_\Delta(m_i)$ and $i=1,\dots,n$, $j\neq i$.
\end{prop}

\begin{proof}

It is obvious that for all $x\in B_\Delta(m_i) $
\be
\label{eq:ww}
\sigma^i(\e)\leq T^i(\e)\leq\tau^i(\e)\quad \Px\text{-a.s.}
\ee
On the other hand, the main contribution to $\tau(\e)$ is made by the switching time
$T(\e)$, for if the trajectory overcomes the saddle point and is in $\Omega_\e^j$ for some $j\neq i$, 
it follows the deterministic trajectory with high probability and reaches the set $B_\Delta(m_j)$ in short
(logarithmic) time.

First we show that
\be
\label{eq:k}
\lim_{\e\downarrow 0}\Px\left(\tau^i(\e)\leq T^i(\e)+ c|\!\ln\e|\right)= 1,
\ee
where $c$ is defined in \eqref{eq:c}.
Let $X_{T^i(\e)}^\e(x)\in\Omega_\e^j$ for some $j\neq i$. On the event 
$A_\e=\{\omega: \sup_{t\in[0,\mu|\!\ln\e|]}|\e L_{t+T^i(\e)}- \e L_{T^i(\e)} |\leq \e^{4\g}\}$ 
the trajectory $X^\e_t(X^\e_{T^i(\e)}(x))$ follows the deterministic trajectory $x^0_t(X^\e_{T^i(\e)}(x))$
which reaches the small neighbourhood of the local minimum $m_j$ in time $c|\!\ln\e|$.
The limit \eqref{eq:k} holds since $\Px(A_\e)\to 1$. Then
% \ba
% \Px&(X_{\tau^i(\e)}^\e\in B_\Delta(m_j)) &\\
% &\leq  \Px(X_{\tau^i(\e)}^\e\in  B_\Delta(m_j),   X_{T^i(\e)}^\e\in \Omega^j_\e,A_\e)
% + \Px(X_{\tau^i(\e)}^\e\in  B_\Delta(m_j),   X_{T^i(\e)}^\e\notin \Omega^j_\e,A_\e)\\
% &\leq \Px( X_{T^i(\e)}^\e\in \Omega^j_\e)
% +\Px(A_\e^c)  \to  \frac{q_{ij}}{q_i}
% \ea
% On the other hand,
\ba
\Px&(X_{\tau^i(\e)}^\e\in B_\Delta(m_j)) 
\geq  \Px(X_{\tau^i(\e)}^\e\in  B_\Delta(m_j),   X_{T^i(\e)}^\e\in \Omega^j_\e,A_\e)\\
&= \Px(  X_{T^i(\e)}^\e\in \Omega^j_\e,A_\e)
\geq \Px( X_{T^i(\e)}^\e\in \Omega^j_\e)-\Px(A_\e^c)  \to  \frac{q_{ij}}{q_i}
\ea
and \eqref{eq:t1} is proved since $\sum_{j\neq i}\frac{q_{ij}}{q_i}=1$.

Convergence \eqref{eq:t2} follows easily from inequality \eqref{eq:ww}, limits \eqref{eq:T2} and \eqref{eq:k} 
and the fact that $\lambda^i(\e)|\!\ln\e|\to 0$.

To prove \eqref{eq:t3} we repeat the argument of Proposition \ref{p:gT}. Indeed, for any $\delta>0$
there
is  $\e_0>0$ such that for $0<\e\leq \e_0$ 
the following inequalities hold
\ba
&\sup_{x\in \Omega^i_\e}\Px\left(  \sup_{t\in[0,c|\!\ln\e|]}|\e L_t |\leq \e^{4\g}\}  \right)
\leq \delta,\\
&\sup_{x\in \Omega^i_\e}\lambda^i(\e)\E_x T^i(\e)\leq 1+\delta,\\
&\max_{1\leq j\leq n-1}\sup_{x\in B_{2\e^\g}(s_j)}\lambda^i(\e)\E_x S^j(\e)\leq \delta,\\
&\max_{1\leq i\leq n} \lambda^i(\e)c|\!\ln \e|\leq \delta.
\ea
Then it is easy to see that for $0<\e\leq \e_0$
\ba
\lambda^i(\e)\E_x \tau^i(\e)&\leq  \lambda^i(\e)(\E_x T^i(\e)+c|\!\ln \e|)
+ \sum_{k=1}^\infty \left[  (1+\delta+ \lambda^i(\e)c|\!\ln \e|)(k+1)+\delta k   \right] \delta^k\\ 
&\leq
1+ \delta\cdot\Const
\ea
which finishes the proof.
\end{proof}

\section{Metastable behaviour. Proof of Theorem \ref{th:main}\label{s:meta}}

\subsection{Convergence on short time intervals}

\begin{prop}
\label{p:m1}
Let $0<\delta<r$. Then if $x\in \Omega^i$, $i=1,\dots,n$, then for $t>0$
\ba
&X^\e_{t/\e^\delta }(x)\stackrel{\mathcal{D}}{\to} m_i,\quad \e\downarrow 0.
\ea

\end{prop}
\begin{proof}
For some $1\leq i\leq n$, let $x\in \Omega^i$.
We shall prove a stronger result: for any $A>0$ and $0<\Delta<\Delta_0$
\be
\label{eq:m2}
\Px\left(\sup_{s\in[c\e^\delta|\!\ln\e|,A]}|X^\e_{s /\e^\delta}-m_i |\leq \Delta \right)
=\Px\left(\sup_{s\in[c|\!\ln\e|,A /\e^\delta]}|X^\e_s-m_i |\leq \Delta \right)\to 1,\quad
 \e\downarrow 0.
\ee
Indeed, recalling Section~\ref{s:onewell} we choose $\g>0$ and $c>0$ such that 
$|X^\e_{c|\!\ln\e|}(x)-m_i|\leq \Delta/2$ a.s.\
on the event 
$E=\mathcal{E}_{c|\!\ln\e|}\cap\{T_1>c|\!\ln\e|\}$, where 
$\mathcal{E}_{c|\!\ln\e|}=\{\sup_{[0,c|\!\ln\e|]}|\e\xi^\e_t|\leq\e^{4\g}\}$. This gives
\ba
\Px&\left(\sup_{s\in[c|\!\ln\e|,A /\e^\delta]}|X^\e_s-m_i |> \Delta \right)\leq
\sup_{|y-m_i|\leq\Delta/2}\P_{\!y}
\left(\sup_{s\in[0,A /\e^\delta-c|\!\ln\e|]}|X^\e_s-m_i |> \Delta \right)+\P(E^c)\\
&\leq \sup_{|y-m_i|\leq\Delta/2}\P_{\!y}
\left(\sigma_\Delta(\e)< A /\e^\delta -c|\!\ln\e| \right)+\P(E^c)\\
&\leq \sup_{|y-m_i|\leq\Delta/2}\P_{\!y}
\left(\sigma_\Delta(\e)< A /\e^\delta \right)
+\P(\mathcal{E}_{\mu|\!\ln\e|}^c)+\P(T_1\leq c|\!\ln\e|),
\ea
where $\sigma_\Delta(\e)=\inf\{t>0:|X^\e_t-m_i|>\Delta  \}$. On the other hand we know
that for $\lambda_\Delta(\e)=H_-(-\Delta/\e)+H_+(\Delta/\e)$,
\be
\lambda_\Delta(\e)\sigma_\Delta(\e)\stackrel{\mathcal{D}}{\to}\exp(1).
\ee
Since $\lambda_\Delta(\e)/\e^\delta\to 0$ as $\e\downarrow 0$
 we have $\P_{\!y}
\left(\sigma_\Delta(\e)< A /\e^\delta \right)\to 0$, as well as
$\P(\mathcal{E}_{c|\!\ln\e|}^c)\to 0$ and $\P(T_1\leq c|\!\ln\e|)\to 0$ 
in the limit of small $\e$. This
 finishes the proof of \eqref{eq:m2}.  
\end{proof}

\begin{rem}
It is easy to notice in view of Section~\ref{s:onewell} that  the 
convergence in Proposition~\ref{p:m1} 
is uniform in $x$ for $x\in \Omega^i_\e$.
\end{rem}

%%%
%%%
\subsection{Proof of Theorem \ref{th:main}}
%%%
%%% 

\begin{lemma}
\label{l:tomin}
For any $t>0$   and $0<\Delta<\Delta_0$,
\be
\Px\left(X^\e_{t/H(1/\e)}\in \bigcup_{j=1}^n B_\Delta(m_j) \right)\to 1, \quad
\e\downarrow 0, 
\ee
uniformly for $x\in\mathbb{R}$.
\end{lemma}
\begin{proof}
Choose $\rho$ and $\g$ such that Proposition \ref{p:q} and Lemma \ref{l:s} hold for small $\e$.
1.\ Let $|x-s_i|\leq 2\e^\g$ for some $i=1, \dots,n-1$.
We know (see Lemma~\ref{l:s}) that
$S_x^i(\e)\leq J(\e)=\inf\{t>0:|\e L_t-\e L_{t-}|>4\e^\g\}$ a.s.\ and 
$ J(\e)\sim \exp(\frac{1}{H(4/\e^{1-\g})})$ if $1-\g>\rho$. 
We show that
\be
\label{eq:tomin}
\Px\left(X^\e_{2/\e^{r(1-\gamma/2)}}\notin \bigcup_{j=1}^n B_\Delta(m_j)\right)\to 0.
\ee
Indeed, the strong Markov property implies
\ba
&\Px\left(X^\e_{2/\e^{r(1-\gamma/2)}}\notin \bigcup_{j=1}^n B_\Delta(m_j) \right)\\
&\leq\Px\left(X^\e_{2/\e^{r(1-\gamma/2)}}\notin \bigcup_{j=1}^n B_\Delta(m_j), 
S^i(\e)\leq 1/\e^{r(1-\gamma/2)} \right)
+\P\left(J(\e)> 1/\e^{r(1-\gamma/2)} \right)\\
&=
\sum_{k=1}^n
\E_x\left[
\P_{\!X^\e_{S^i(\e)}}\left( X^\e_{2/\e^{r(1-\gamma/2)}-S^i(\e)}
\notin \bigcup_{j=1}^n B_\Delta(m_j)\right)\cdot 
\I{S^i(\e)\leq 1/\e^{r(1-\gamma/2)}}\cdot \I{X^\e_{S^i(\e)}\in \Omega^k_\e}\right] \\
&+\E_x\left[
\P_{\!X^\e_{S^i(\e)}}\left( X^\e_{2/\e^{r(1-\gamma/2)}-S^i(\e)}
\notin \bigcup_{j=1}^n B_\Delta(m_j) \right)\cdot 
\I{S^i(\e)\leq 1/\e^{r(1-\gamma/2)}}\cdot \I{X^\e_{S^i(\e)}\notin \bigcup_{j=1}^n \Omega^j_\e }\right] \\
&+\P\left(J(\e)> 1/\e^{r(1-\gamma/2)} \right)\\
&\leq
\sum_{k=1}^n\sup_{y\in \Omega_\e^k}\P_{\!y}\left(\sup_{s\in [c|\!\ln\e|,2/\e^{r(1-\gamma/2)}]} 
|X^\e_s-m_k|>\Delta\right)\\
&+\P\left(\sup_{t\in[0, 1/\e^{r(1-\gamma/2)}   ]} \e |L_{t}-L_{t-}|> a \right) 
+\P\left(J(\e)> 1/\e^{r(1-\gamma/2)} \right),
\ea 
with $a=\frac{1}{2}\min\{s_2-s_1, \dots, s_{n-1}-s_{n-2}\}$.
The first summand in the latter formula tends to $0$ due to Proposition~\ref{p:m1}. 
The second summand is estimated by $1-\exp(\e^{-r(1-\gamma/2)}H(a/\e))\to 0$, and the third summand also
tends to $0$ due to the definition of $J(\e)$. 

2.\ It is clear from the proof that the limit \eqref{eq:tomin} holds also for $x\in \Omega_\e^i$, $i=1,\dots,n$, 
and thus for all $x\in\mathbb{R}$.
Then, for $\e$ small enough such that $t/H(1/\e)>2/\e^{r(1-\g/2)}$ the application of the Markov property 
\be
\Px\left(X^\e_{t/H(1/\e)}\notin  \bigcup_{j=1}^n B_\Delta(m_j) \right)=
\E_x\P_{\! X^\e_{t/H(1/\e)   -2/\e^{r(1-\g/2)}}}
\left(X^\e_{2/\e^{r(1-\g/2)}}\notin \bigcup_{j=1}^n B_\Delta(m_j) \right)
\ee 
finishes the proof.
\end{proof}

\bigskip

\begin{Aproof}{Theorem}{\ref{th:main}}

It is clear from the Markov property that it is sufficient to show that for any $t>0$
and $x\in \Omega^i_\e $, $i=1,\dots, n$,
\ba
\Px\left(X^\e_{t/H(1/\e)}\in B_\Delta(m_j)\right)\to \P_{\!m_i}\left( Y_t=m_j\right),
\quad j\neq i.
\ea

Define a sequence of stopping times $(\tau(k))_{k\geq 0}$ and states $(m(k))_{k\geq 0}$ such that
$\tau(0)=0$, $m(0)=m_i$ and for $k\geq 1$
\ba
\tau(k)&=\inf\{t>\tau(k-1):X^\e_t\in \bigcup_{i=1}^n B_\Delta(m_i)\backslash B_\Delta(m(k-1))\},\\
m(k)&=\sum_{j=1}^n m_j\I{X^\e_{\tau(k)}\in B_\Delta(m_j)} .
\ea

Define also a (non-Markovian) process $\hat X^\e$ on a state space $\{m_1,\dots,m_n\}$
\be
\hat X^\e_t=
\sum_{k=0}^\infty    m(k)\cdot\I{t\in\left[ H(1/\e)\tau(k),H(1/\e)\tau(k+1)\right) }.
\ee

The strong Markov property of $X^\e$ and Proposition~\ref{p:1f} imply that 
\ba
\label{eq:12}
&\Law\left(H(1/\e)(\tau(k+1)-\tau(k)) |\hat{X}^\e_{\tau(k)}=m_i\right) \to  \exp(1/q_i),\\
&\Px(\hat{X}^\e_{\tau(k+1)}=m_j|\hat{X}^\e_{\tau(k)} =m_i)\to \frac{q_{ij}}{q_i},
\ea
uniformly for $k\geq 0$.

The process $Y$ defined in the statement of the Theorem is given by the sequence of its jump times and states,
$(\theta(k), Y_k )_{k\geq 0}$ with the property that
the interjump times are conditionally independent and exponentially distributed.
% Denote
% \be
% \lambda_i=\lim_{\e\downarrow 0}\frac{\lambda_i(\e)}{H(1/\e)}=|s_{i-1}-m_i|^{-r}+|s_i-m_i|^{-r}.
% \ee
% Then 
and for $k\geq 0$, $0\leq i,j\leq n$, $i\neq j$,
\ba
\label{eq:13}
&\Law\left(\theta(k+1)-\theta(k)|Y_{k}=m_i\right)  =\exp(1/q_i)\\
&\P(Y_{k+1}=m_j|Y_k=m_i)= \frac{q_{ij}}{q_i}.
\ea
Then %for $\e$ small enough such that  $x\in \Omega^i_\e$
\ba
\label{eq:11}
&\Big|\Px\left(X^\e_{t /H(1/\e)}\in B_\Delta(m_j)\right)- 
\P_{\!m_i}\left(Y_t=m_j \right)\Big|\\
&\leq 
\Big| \Px\left(  X^\e_{t /H(1/\e)}\in B_\Delta(m_j)      \right)- 
\Px\left(   \hat{X}^\e_t=m_j        \right)\Big|
+\Big|\Px\left(\hat{X}^\e_t=m_j \right)-\P_{\!m_i}\left(Y_t=m_j \right)\Big|.
\ea
The second summand in \eqref{eq:11} vanishes in the limit of small $\e$ due to the weak
convergence of the jump process $\hat{X}^\e$ to $Y$. Indeed, in this case the weak convergence
is equivalent to the weak convergence of the sequences of jump times and jump sizes (see \cite{Xia-92})
$(\tau(k),m(k))_{k\geq 0}\Rightarrow (\theta(k),Y_k)_{k\geq 0}$,
which follows from \eqref{eq:12} and \eqref{eq:13}.

To estimate the first summand in \eqref{eq:11} we use Lemma~\ref{l:tomin}. Indeed,
\ba
&\Big| \Px\left(  X^\e_{t /H(1/\e)}\in B_\Delta(m_j)      \right)- 
\Px\left(   \hat{X}^\e_t=m_j        \right)\Big|
\\
&=\Big|\Px\left(   X^\e_{t /H(1/\e)}\in B_\Delta(m_j)  , \hat{X}^\e_t=m_j\right)
+\Px\left(   X^\e_{t /H(1/\e)}\in B_\Delta(m_j)       , \hat{X}^\e_t\neq m_j\right)(=0)\\
&-\Px\left(\hat{X}^\e_t=m_j,   X^\e_{t /H(1/\e)}\in  B_\Delta(m_j)        \right)
-\Px\left(\hat{X}^\e_t=m_j, X^\e_{t /H(1/\e)}\in \bigcup_{k\neq j} B_\Delta(m_k) \right)(=0)\\
&-\Px\left(\hat{X}^\e_t=m_j, X^\e_{t /H(1/\e)}\notin \bigcup_{k= 1}^n B_\Delta(m_k) \right)\Big|
\leq \Px\left(   X^\e_{t /H(1/\e)}\notin \bigcup_{k=1}^n B_\Delta(m_k)     \right)\to 0,
\ea
which finishes the proof of the Theorem.
\end{Aproof}

\begin{appendix}
\section{Existence of strong solution\label{a:ss}}

%\appendix

Here we refer to
\cite{SamorodnitskyG-03} where the existence of the strong solution was established for 
potentials with unique stable point.

First, we note that for the existence and uniqueness of solutions of \eqref{eq:main} it is
enough to demand that $U'$ is locally Lipschitz and $U'(x)x\geq 0$ for $|x|\geq N$, with
$N$ large enough.

For brevity, we set $\e=1$. Then for $n\geq 1$ consider a family of SDEs with truncated drift,
\be
\label{eq:mainc}
X^{(n)}_t(x)=x-\int_0^t U'([X^{(n)}_{s-}(x)]_n)\, ds+L_t,
\ee
where 
\be
[x]_n=
\begin{cases}
&-n^2,\quad x< -n^2\\
&\hphantom{-}x,\quad -n^2\leq x\leq n^2,\\
&\hphantom{-}n^2,\quad  x>n^2.
\end{cases}
\ee
Since for any $n\geq 1$ the drift $U'([\,\cdot\,]_n)$ is globally Lipschitz,
\eqref{eq:mainc} has a strongly unique
solution for any $\mathcal{F}_0$-measurable $x$ which is a semimartingale 
by \cite[Theorem V.3.6]{Protter-04} and also strongly
Markov by \cite[Theorem V.6.34]{Protter-04}.

%In order to obtain a strong and unique solution of \eqref{eq:main} for any $x\in\mathbb{R}$ 
%we show according
%to \cite[Theorem V.7.34]{Protter-04} that the explosion time $T_x(\omega)=\infty$ a.s. 

Define a family of stopping times $T_x^n:\Omega\to\mathbb{R}\cup\{+\infty\}$
in the following way
\be
T_x^{n}=\inf\{t>0: |X^{(n)}_t(x)|> n^2\},  %  \text{ or } |X^{(n)}_{t-}(x)|\geq n^2\}, 
\quad n\geq 1,
\ee
with the usual convention that if $|X^{(n)}_t(x,\omega)|\leq n^2$ for all $t\geq 0$, then 
$T^n_x(\omega)=+\infty$.
As $(T^n_x)_{n\geq 1}$ is a non-decreasing sequence, we can define for any $\omega\in\Omega$
the \textit{explosion time}
\be
T_x(\omega)=\lim_{n\to\infty}T_x^n(\omega).
\ee
If we show that $T_x=+\infty$ a.s.\ then using the fact that $X^{(n)}_t(x)=X^{(n+1)}_t(x)$
a.s.\ for $0\leq t<T^n_x$ we obtain the solution of \eqref{eq:main} by setting
\be
X_t(x)=X^{(n)}_{t\wedge T_x^n}(x)
\ee
for all $t$ and $\omega$ such that $t<T_x^n$, see \cite[Theorem V.7.34]{Protter-04}.

Let us suppose to the contrary that there is $0<A<\infty$ such that for some $x\in \mathbb{R}$
\be
\P\left(\lim_{n\to\infty} T_x^n\leq A\right)=\delta>0.
\ee
According to this choose $B>0$ such that
\be
\P\left(\sup_{t\in[0,A]}L_t\leq B\right)>1-\delta.
\ee
Fix $\omega\in\{   \lim_{n\to\infty} T_x^n\leq A \}
\cap \{   \sup_{t\in[0,A]}L_t\leq B \}$.
For $n\geq 1$ let
%\be
%S_x^n=\sup\{t\in[T_x^n,T_x^{n+1}):\sgn X^{n+1}_t(x)=  \sgn X^{n+1}_{T^{n+1}_x}(x) \text{ and }
%|X^{n+1}_t(x)|\geq N \text{ on } [S_x^n, T_x^{n+1}) \}
%\ee
\be
S_x^n=\inf\{t\in[T_x^n,T_x^{n+1}):\sgn X^{n+1}_s(x)=  \sgn X^{n+1}_{T^{n+1}_x}(x)  \text{ and }
|X^{n+1}_s(x)|\geq N \text{ for } s\in[t, T_x^{n+1}) \}
\ee
and $S_x^n=T_x^{n+1}-$, if there is $t_0\in [T_x^n,T_x^{n+1})$ such that for 
$s\in [t_0, T_x^{n+1})$
we have $\sgn X^{n+1}_s(x)\neq  \sgn X^{n+1}_{T^{n+1}_x}(x)$ or
$|X^{n+1}_s(x)|< N$ (or both).

Consider different cases separately. First, if $S_x^n=T_x^{n+1}-$, we have necessarily that the 
jump size  $|L_{T^{n+1}_x}-L_{T^{n+1}_x-}|$ must be bigger than
$(n+1)^2-N$ which leads to contradiction for $n$ big enough. If $S_x^n\neq T^{n+1}_x$ we have that
$|X^{n+1}_{S^n_x}(x)|\leq n^2+2B$ and 
$|X^{n+1}_{T^{n+1}_x}(x)|\geq (n+1)^2$. 
%which also leads to contradiction for $n$ big enough. 
%Finally, if $S^n_x>T^{n}_x$ 
We note that due to the inequality $U'(x)x\geq 0$ for $|x|\geq N$, 
we have
\be
\sgn\left(  X^{(n+1)}_{T^{n+1}_x}(x)-X^{(n+1)}_{S^n_x}(x)\right) =
\sgn \int_{S^n_x}^{T^{n+1}_x}U'([X^{(n+1)}_{s-}(x) ]_n)\,ds 
\ee
and
\ba
|L_{T^{n+1}_x}-L_{S^n_x}|&=
\left| X^{(n+1)}_{T^{n+1}_x}(x)-X^{(n+1)}_{S^n_x}(x)+
\int_{S^n_x}^{T^{n+1}_x}U'([X^{(n+1)}_{s-}(x) ]_n)\,ds \right|\\
&\geq 
|X^{(n+1)}_{T^{n+1}_x}(x)-X^{(n+1)}_{S^n_x}(x)|
\geq (n+1)^2-n^2-2B,
\ea
which also contradicts the assumptions for $n$ sufficiently big. Thus the existence and uniqueness of the strong 
solution of \eqref{eq:main} is established. This solution is also strongly Markov and Feller, see 
\cite{SamorodnitskyG-03} and \cite[Theorem 5.4]{KurtzP-91}

\section{Regular variation\label{a:rv}}

\begin{defin}
a) A positive, Lebesgue measurable function $l$ on $(0,+\infty)$ is slowly varying at $+\infty$ if
\be
\label{eq:sv}
\lim_{u\to+\infty}\frac{l(\lambda u)}{l(u)}=1,\quad \lambda>0.
\ee
b) A positive, Lebesgue measurable function $H$ on $(0,+\infty)$ is regularly varying at $+\infty$ of index 
$r\in\mathbb{R}$ if
\be
\lim_{u\to+\infty}\frac{H(\lambda u)}{H(u)}=\lambda^r,\quad \lambda>0.
\ee
\end{defin}
For example, positive constants, logarithms and iterated logarithms are slowly varying functions.
Further, one can prove that $H$ is regularly varying of index $r$ if and only if there is 
a slowly varying function $l$ such that
\be
H(u)=u^r l(u), \quad u>0.
\ee
Another important result is the uniform convergence in \eqref{eq:sv}.
\begin{prop}[\cite{BinghamGT-87}, Theorem 1.2.1]
\label{p:usv}
If $l$ is slowly varying at $+\infty$ then 
\be
\lim_{u\to+\infty} \frac{l(\lambda u)}{l(u)}=1,
\ee  
uniformly for $\lambda$ from a compact set in $(0,+\infty)$. 
\end{prop}

\end{appendix}

\bibliography{biblio-new}
\bibliographystyle{alpha}

\vspace{1cm}

\parbox{.48\linewidth}{
\noindent
Peter Imkeller\\
Institut f\"ur Mathematik\\
Humboldt Universit\"at zu Berlin\\
Rudower Chaussee 25\\
12489 Berlin Germany\\
E.mail: imkeller@mathematik.hu-berlin.de\\
http://www.mathematik.hu-berlin.de/\~{}imkeller
}\hfill
\parbox{.48\linewidth}{
\noindent
Ilya Pavlyukevich\\
Institut f\"ur Mathematik\\
Humboldt Universit\"at zu Berlin\\
Rudower Chaussee 25\\
12489 Berlin Germany\\
E.mail: pavljuke@mathematik.hu-berlin.de\\
http://www.mathematik.hu-berlin.de/\~{}pavljuke
}

\end{document}